\documentclass[11pt,reqno]{amsart}
\usepackage[a4paper,textwidth=158mm,textheight=235mm,centering]{geometry}
\usepackage{amsmath,amssymb,amsthm,mathtools,bm}
\usepackage[T1]{fontenc}
\usepackage{lmodern}
\usepackage[unicode,hidelinks]{hyperref}
\numberwithin{equation}{section}
\newtheorem{theorem}{Theorem}[section]
\newtheorem{proposition}[theorem]{Proposition}
\newtheorem{lemma}[theorem]{Lemma}
\newtheorem{corollary}[theorem]{Corollary}
\theoremstyle{remark}
\newtheorem{remark}[theorem]{Remark}
\newcommand{\R}{\mathbb R}
\newcommand{\HH}{\mathbf H}
\newcommand{\Res}{\mathcal R}
\newcommand{\dd}{\,d}
\newcommand{\tr}{\operatorname{tr}}
\newcommand{\Div}{\operatorname{div}}
\newcommand{\spt}{\operatorname{spt}}
\newcommand{\dist}{\operatorname{dist}}
\newcommand{\Fcal}{\mathcal F}
\allowdisplaybreaks[2]
\title[Bounded total curvature II: multiplicity one]{Singularities of Mean Curvature Flow with Bounded Total Curvature II: Multiplicity One}
	\author{Yongheng Han}
	\address{Institute of Mathematics, Academy of Mathematics and Systems Science, Chinese Academy of Sciences, Beijing 100190, China}

	\email{hyh98@amss.ac.cn}

    \author{Zhen Wang}
	\address{The Institute of Geometry and Physics, University of Science and Technology
of China, Hefei,  Anhui 230026, China}

	\email{wz1854@ustc.edu.cn}
\date{\today}

\begin{document}

\begin{abstract}
Let \(M_t^n\subset\mathbb R^{n+1}\) be a smooth mean curvature flow of closed embedded hypersurfaces, with \(n\ge3\), and let \(T\) be its first singular time. If the integral of \(|A|^n\) is uniformly bounded for \(t<T\), every nonzero tangent flow at time \(T\) has multiplicity one. At each singular point the tangent flow is unique, and the rescaled flow converges locally smoothly to a closed self-shrinker. 
\end{abstract}
\maketitle

\section{Introduction}

Let \(F:M^n\times[0,T)\to\R^{n+1}\) be a smooth mean curvature flow of closed embedded hypersurfaces, and let \(M_t=F(M,t)\). Our sign convention is
\begin{equation}\label{eq:intro-mcf}
  \partial_tF=\HH.
\end{equation}
 Assume that \(T<\infty\) is the first singular time and that
\begin{equation}\label{eq:total}
  \sup_{0\le t<T}\int_{M_t}|A|^n\dd\mu_t\le\Lambda<\infty .
\end{equation}
Here \(A\) denotes the second fundamental form and \(\mu_t\) the induced area measure.
For \(x_0\in\R^{n+1}\) and \(\lambda>0\), set
\begin{equation}\label{eq:intro-rescaling}
  M^\lambda_\tau=\lambda(M_{T+\lambda^{-2}\tau}-x_0),
  \qquad -\lambda^2T<\tau<0 .
\end{equation}
A tangent flow at \((x_0,T)\) is a subsequential Brakke limit of these flows as \(\lambda\to\infty\). Throughout this paper tangent flows are considered at negative times. By Huisken's monotonicity formula, every nonzero tangent flow is backward self-similar \cite{Brakke,Huisken,IlmanenCompact}. On the regular part of its time \(-1\) slice, the support \(\Sigma\) satisfies
\begin{equation}\label{eq:intro-shrinker}
   \HH_\Sigma+\tfrac12 x^\perp=0 .
\end{equation}

Ilmanen proved smoothness of tangent-flow supports at the first singular time for closed embedded surfaces in \(\R^3\) and formulated the multiplicity-one conjecture \cite{Ilmanen95}. White proved that tangent flows at the first singular time of a mean-convex flow have multiplicity one \cite[Theorem~1]{White03}. For embedded surface flows, Li and Wang proved the blow-up of the mean curvature at the first singular time and established multiplicity one under a type-I bound on the mean curvature \cite{LW19,LW22}. Sun studied perturbations that exclude a prescribed closed singularity model of higher multiplicity \cite{Sun23}. Bamler and Kleiner proved the multiplicity-one conjecture for surfaces without a curvature assumption \cite{BK}.

For higher-dimensional hypersurfaces satisfying \eqref{eq:total}, we proved that the support of every tangent flow at the first singular time is smooth \cite[Theorems~1.1 and~5.1]{HW}. The integer multiplicity on this support is the subject of the present paper.

\begin{theorem}\label{thm:main}
Let \(n\ge3\). For a smooth mean curvature flow of closed embedded hypersurfaces in \(\R^{n+1}\) satisfying \eqref{eq:total}, every nonzero tangent flow at the first singular time has multiplicity one.
\end{theorem}

\begin{corollary}\label{cor:unique}
Under the assumptions of Theorem~\ref{thm:main}, the tangent flow at each singular point \((x_0,T)\) is unique. Its time \(-1\) slice is a smooth closed embedded self-shrinker. As \(t\uparrow T\), the hypersurfaces \((T-t)^{-1/2}(M_t-x_0)\) converge locally smoothly to this self-shrinker with multiplicity one.
\end{corollary}

Theorem~\ref{thm:main} and Corollary~\ref{cor:unique} address \cite[Conjecture~1.3]{HW} in the closed embedded setting. We isolate the step that uses the separation of sheets.

\begin{proposition}\label{prop:conditional}
Let \(n\ge3\), and suppose a smooth closed embedded flow satisfies \eqref{eq:total}. Fix \((x_0,T)\) admitting a nonzero tangent flow. If the support of every nonzero tangent flow there is smoothly and properly embedded at each negative time, then every such flow has multiplicity one.
\end{proposition}

Put \(s=-\log(T-t)\). Along a subsequence of translates of the normalized flow, the limit on a fixed unit time interval is a stationary rescaled self-shrinker. The curvature-varifold compactness theorem recalled in \cite[Lemma~2.13]{HW} passes the critical curvature bound to this limit. On each tubular chart the integrated Bochner formula for the signed distance also gives strong local spacetime \(L^2\) convergence of the second fundamental forms. In particular, curvature is measured against the area of all the sheets in the limit. We give the details in Section~\ref{sec:limits}.

For \(x\in M_t\), let \(\mathfrak s(x,t)\) be the infimum of those \(r>0\) for which \(M_t\cap B_r(x)\) is disconnected, with value \(+\infty\) if no such radius exists. This is the separation function of Bamler and Kleiner \cite[Definition~3.1]{BK}. Here and throughout, \(B_r(x)\) is a Euclidean ball in the ambient space. Suppose that a sequence of normalized slices converges to \(m|\Sigma|\) with \(m\ge2\). Away from the atoms of a fixed positive size in the limiting curvature measure, the component through the center of a small ball has area close to \(\omega_n r^n\), by \cite[Corollary~4.4]{HW}. The entire intersection has limiting area close to \(m\omega_n r^n\). Thus the ball contains another component. Taking \(r\) arbitrarily small shows that the normalized separation tends to zero in Gaussian measure.

The evolution of the logarithmic separation was studied in \cite[Section~3]{BK}. On a region where the component through the base point is a graphical disk and the separation is realized by the distance to its complement, the two-point calculation gives
\[
 (\partial_s-\Delta_{M_s})\log\widehat{\mathfrak s}\ge\tfrac12,
 \qquad \widehat{\mathfrak s}=(T-t)^{-1/2}\mathfrak s .
\]
This local statement requires care near points where the connected components of a ball intersection change. We introduce a spacetime length \(\ell\) defined by connectedness and curvature control on parabolic cylinders. In regions modeled on a sheet of a collection of parallel planes, \(\ell\) approaches the distance to the nearest other plane. Near the remaining points we replace the separation by a smooth positive length \(b\) comparable to \(\ell\). The construction yields a single function on the flow for every fixed \(0<\delta<1\),
\[
 h_\delta=\log\min\{\delta,\widehat{\mathfrak s},\widehat b\},
 \qquad \widehat b=(T-t)^{-1/2}b.
\]

The use of a smoothed local scale to modify the separation, followed by integration against a backward heat kernel, follows the strategy of \cite[Sections~4--6]{BK}. To estimate the cost of this modification in higher dimensions, we use the bound \eqref{eq:total}. In ordinary flow coordinates, on a cylinder of spatial scale \(r\), the contribution of \(\ell^{-2}\) is of order at most \(r^n\) before insertion of the Gaussian factor. A compactness argument shows that each relevant cylinder either contributes to a drop of Gaussian area or satisfies
\[
             \int_Q |A|^n\dd\mu_t\,dt>r^2
\]
on a fixed enlargement \(Q\). In the latter case the cost is bounded by \(r^{n-2}\int_Q |A|^n\dd\mu_t\,dt\). At any given scale the enlarged cylinders have bounded overlap. The weights at smaller scales form a convergent geometric series because \(n>2\). At the cylinders controlled by Gaussian area, bounded overlap across scales follows from monotonicity centered at the points of integration. These two estimates allow the negative part in the evolution of \(h_\delta\) to be made small on every sufficiently late normalized unit interval.

If one tangent flow at the fixed point were multiple, the Gaussian integral of \(h_\delta\) would then increase by a definite positive amount on every such interval for one fixed \(\delta\). The upper bound \(h_\delta\le\log\delta\) rules this out. This argument estimates the region where the local graphical description can fail by integration in spacetime.

Theorem~\ref{thm:main} also admits a proof using the two one-sided
tangent-ball radii. Choose a smooth unit normal \(\nu\) along the
rescaled flow and set \(\nu_+=\nu\) and \(\nu_-=-\nu\). For
\(x\in M_s\) and \(\alpha\in\{+,-\}\), define
\[
 r_\alpha(x,s)
   =\sup\bigl\{a>0:
       B_a\bigl(x+a\nu_\alpha(x,s)\bigr)\cap M_s=\varnothing\bigr\}.
\]
Thus \(r_\alpha(x,s)\) is the supremal radius of an open ball tangent
to \(M_s\) at \(x\) on the side \(\nu_\alpha\) and disjoint from
\(M_s\). Fix a sufficiently large dimensional constant \(K\). For
\(0<\delta<1\), define
\[
 q_\delta(x,s)
   =\delta\exp\!\left[-K\delta^{\,n-2}
          \bigl(1+|A(x,s)|^2\bigr)^{(n-2)/2}\right],
 \qquad
 h_{\delta,\alpha}(x,s)
   =\log\min\{r_\alpha(x,s),q_\delta(x,s)\}.
\]
The choice of \(K\) ensures that \(q_\delta|A|<1\). Hence, wherever
\(r_\alpha<q_\delta\), the tangent ball is limited by a second point
of contact with \(M_s\). The evolution of its radius can then be
estimated by a two-point calculation. If a tangent flow has
multiplicity greater than one, the nearby sheets force at least one
of the two radii below \(q_\delta\) on a set of positive Gaussian
measure. Summing the estimates for the two sides and using the total
curvature bound to control the error from \(q_\delta\), one finds
that, for a fixed sufficiently small \(\delta\), the quantity
\[
 J_\delta(s)
   =\int_{M_s}\bigl(h_{\delta,+}+h_{\delta,-}\bigr)
        (4\pi)^{-n/2}e^{-|x|^2/4}\,\dd\mu_s
\]
increases by a fixed positive amount on every sufficiently late unit
time interval. This contradicts \(J_\delta(s)\le0\), which follows
from \(q_\delta\le\delta<1\).

\subsection*{Organization of the paper}

In Section~\ref{sec:limits} we establish the curvature convergence, describe the possible limiting self-shrinkers, and use the componentwise area estimate to detect multiple sheets in ambient balls. In Section~\ref{sec:contact} we study the local scale and its behavior at planar limits. We then derive the two-point evolution inequality and construct the smooth truncation used to integrate it. In Section~\ref{sec:proof} we estimate the Gaussian integral on the region where the separation has been modified. We prove Proposition~\ref{prop:conditional} and deduce Theorem~\ref{thm:main}, uniqueness of the tangent flow, and local smooth convergence at each singular point.

\subsection*{Disclosure of AI use}
This paper was generated using ChatGPT-6 Astra. OpenAI Codex (GPT-6) was used in developing and writing the separation-function argument, working out the estimates, and checking the calculations and references for this version. The authors take full responsibility for the content of this article, including its mathematical claims, arguments, and conclusions.

\subsection*{Acknowledgments}
The authors thank Haozhao Li for helpful discussions. This work was supported by the Project of Stable Support for Youth Team in Basic Research Field, Chinese Academy of Sciences (YSBR-001).

\section{Curvature convergence and tangent flows}\label{sec:limits}

\subsection{Notation and conventions.}\leavevmode\par\nobreak

We use the conventions
\begin{equation}\label{eq:sign}
 S(u)=D_u\nu,\qquad B(u,v)=-\langle Su,v\rangle\nu,
 \qquad H=\tr S,\qquad \HH=-H\nu.
\end{equation}
Both arguments of $B$ are extended to the ambient space by orthogonal
projection onto the tangent space. We likewise extend $S$ by zero on the
normal line. Thus $S=-B\cdot\nu$, $|B|=|A|$, and all tensor differences
below are differences of ambient tensors. Let $P$ and $P^\perp$ denote
the tangent and normal projections, respectively.

Let $E$ denote the entropy of the initial hypersurface,
\[
 E:=\sup_{y\in\R^{n+1},\,a>0}
 (4\pi a)^{-n/2}\int_{M_0}e^{-|x-y|^2/(4a)}\dd\mu_0.
\]
Since $M_0$ is smooth and closed, $E$ is finite.
By Huisken's monotonicity formula,
\begin{align*}
 &(4\pi a)^{-n/2}\int_{M_t}e^{-|x-y|^2/(4a)}\dd\mu_t\\
 &\hspace{12mm}\le
 (4\pi(t+a))^{-n/2}\int_{M_0}e^{-|x-y|^2/(4(t+a))}\dd\mu_0.
\end{align*}
Taking the supremum gives
\begin{equation}\label{eq:entropy}
 \sup_{0\le t<T}\sup_{y\in\R^{n+1},\,a>0}
 (4\pi a)^{-n/2}\int_{M_t}e^{-|x-y|^2/(4a)}\dd\mu_t\le E.
\end{equation}
Since $e^{-|x-y|^2/(4r^2)}\ge e^{-1/4}$ on $B_r(y)$, we obtain
\begin{equation}\label{eq:area}
 \mu_t(B_r(y))\le Dr^n,
 \qquad D=e^{1/4}(4\pi)^{n/2}E.
\end{equation}
Under dilation by $\lambda$, the quantities $|A|$, $\dd\mu$, and $r$
become $\lambda^{-1}|A|$, $\lambda^n\dd\mu$, and $\lambda r$.
Thus \eqref{eq:total} and \eqref{eq:area} hold with the same constants
after parabolic rescaling.

Fix a point \((x_0,T)\) such that there exists a nonzero tangent flow. Set
\begin{equation}\label{eq:normalization}
 s=-\log(T-t),\qquad X=e^{s/2}(F-x_0),\qquad
 \rho(X)=(4\pi)^{-n/2}e^{-|X|^2/4},\qquad
 \Res=\HH+\tfrac12X^\perp.
\end{equation}
For $R\ge1$, the annuli
$B_{k+1}\setminus B_k$, $k\ge\lfloor R\rfloor$, and
\eqref{eq:area} give
\[
 \int_{M_s\setminus B_R}\rho\dd\mu_s
 \le(4\pi)^{-n/2}D
       \sum_{k\ge\lfloor R\rfloor}(k+1)^ne^{-k^2/4}
 \longrightarrow0
 \quad\hbox{as }R\to\infty,
\]
uniformly in normalized time.

Since $dt/ds=e^{-s}$ and the mean curvature vector scales by
$e^{-s/2}$, the normalized hypersurfaces, again denoted by $M_s$, satisfy
\begin{equation}\label{eq:velocity}
 \partial_sX=\HH+\tfrac12X,\qquad \Box=\partial_s-\Delta.
\end{equation}
Here $\partial_s$ denotes differentiation at fixed $p$ in $X(p,s)$.
For the scalar calculation below choose a unit normal locally, with the sign convention in \eqref{eq:sign}. The evolution of the metric and the vector-valued tensors is independent of this choice.
Differentiating $g_{ij}=\langle\partial_iX,\partial_jX\rangle$, we have
\[
 \partial_sg_{ij}=-2HS_{ij}+g_{ij},\qquad
 \partial_s\dd\mu_s
 =\tfrac12g^{ij}\partial_sg_{ij}\dd\mu_s
 =(n/2-H^2)\dd\mu_s.
\]
Moreover,
\[
 \partial_s\rho=-\tfrac12\langle X,\HH+X/2\rangle\rho,
 \qquad \nabla\rho=-\tfrac12X^T\rho.
\]
The identity $\Delta X=\HH$ yields
\[
 \Delta\frac{|X|^2}{2}=n+\langle X,\HH\rangle,
 \qquad \Div_{M_s}X^T=n+\langle X,\HH\rangle.
\]
Hence
\begin{align*}
 \Div_{M_s}(\rho X^T/2)
 &=\left(\frac n2+\frac12\langle X,\HH\rangle
                    -\frac14|X^T|^2\right)\rho,\\
 |\Res|^2
 &=|\HH|^2+\langle X,\HH\rangle+\frac14|X^\perp|^2.
\end{align*}
Subtracting the second expression from the first gives
\begin{align}\label{eq:weighted-volume}
 \partial_s(\rho\dd\mu_s)
 &=\left(\frac n2-|\HH|^2-\frac12\langle X,\HH\rangle
                         -\frac14|X|^2\right)\rho\dd\mu_s\notag\\
 &=\left[-|\Res|^2\rho+
              \Div_{M_s}(\rho X^T/2)\right]\dd\mu_s.
\end{align}
For a smooth function $f$ on the evolving hypersurface, integration by
parts therefore gives
\begin{align}\label{eq:weighted-test}
 \frac d{ds}\int_{M_s}f\rho\dd\mu_s
 &=\int_{M_s}\left(\partial_sf
       -\tfrac12\langle X^T,\nabla f\rangle-f|\Res|^2\right)
            \rho\dd\mu_s\notag\\
 &=\int_{M_s}(\Box f-f|\Res|^2)\rho\dd\mu_s.
\end{align}
For the last equality we used
\[
 \int_{M_s}\Delta f\rho\dd\mu_s
 =-\int_{M_s}\langle\nabla f,\nabla\rho\rangle\dd\mu_s
 =\frac12\int_{M_s}\langle\nabla f,X^T\rangle\rho\dd\mu_s.
\]
If $\phi$ is a fixed smooth ambient function, then
$\partial_s(\phi\circ X)=\bar\nabla\phi\cdot(\HH+X/2)$ and
$\nabla^{M_s}\phi=P\bar\nabla\phi$. Substitution into the first line of
\eqref{eq:weighted-test} gives
\begin{equation}\label{eq:local-mass}
 \frac d{ds}\int_{M_s}\phi\rho\dd\mu_s
 =\int_{M_s}\bar\nabla\phi\cdot\Res\,\rho\dd\mu_s
       -\int_{M_s}\phi|\Res|^2\rho\dd\mu_s.
\end{equation}
In particular,
\begin{equation}\label{eq:F}
 \Fcal(s):=\int_{M_s}\rho\dd\mu_s\downarrow\Theta>0,
 \qquad \Fcal'(s)=-\int_{M_s}|\Res|^2\rho\dd\mu_s.
\end{equation}
The limit is positive because a nonzero tangent flow has positive mass
on some fixed ball, where $\rho$ has a positive lower bound. It follows that
\begin{equation}\label{eq:defect}
 \int_s^{s+1}\int_{M_u}|\Res|^2\rho\dd\mu_u\,du
 =\Fcal(s)-\Fcal(s+1)\longrightarrow0.
\end{equation}

\subsection{Lower semicontinuity and strong curvature convergence}\leavevmode\par\nobreak

\begin{lemma}\label{lem:curvature-lsc}
Let $M_j$ be smooth, properly embedded hypersurfaces without relative boundary in an open
set $U$, and suppose that $|M_j|\rightharpoonup\theta|\Sigma|$ locally,
where $\Sigma$ is smooth and embedded and $\theta$ is a locally constant
positive integer on $\Sigma$. Suppose that $p>1$ and
$\sup_j\int_{M_j\cap K}|A_j|^p\dd\mu_j<\infty$ for every compact
$K\Subset U$. Then every nonnegative $\zeta\in C_c(U)$ satisfies
\begin{equation}\label{eq:lsc}
 \int_\Sigma\zeta\theta|A_\Sigma|^p\dd\mu_\Sigma
 \le\liminf_j\int_{M_j}\zeta|A_j|^p\dd\mu_j.
\end{equation}
\end{lemma}

\begin{proof}
Apply the curvature-varifold compactness theorem of Hutchinson and its extension by Mantegazza \cite{Hutchinson,Mantegazza}, in the formulation of \cite[Lemma~2.13]{HW}. Its boundary measure vanishes here,
because the approximating hypersurfaces have no relative boundary in
$U$. Local mass bounds follow from varifold convergence. To verify the curvature hypothesis and identify the norm in the conclusion, recall that
\[
 \mathcal A_{ijk}=D_{Pe_i}P_{jk},\qquad
 P=I-\nu\otimes\nu.
\]
For a smooth hypersurface,
\[
 D_{Pe_i}P=-(S Pe_i)\otimes\nu-\nu\otimes(S Pe_i),
 \qquad |\mathcal A|^2=2|A|^2.
\]
Thus the entrywise tensor norm used in the hypothesis of that theorem
is locally bounded in $L^p$ by finite-dimensional norm equivalence.
On the smooth limit with locally constant multiplicity, the generalized
curvature tensor is the classical derivative of its tangent projection.
Use the convex integrand
$2^{-p/2}|\mathcal A|^p$ in the lower-semicontinuity conclusion.
Localization, followed by approximation of the nonnegative continuous
weight $\zeta$, gives \eqref{eq:lsc}. In particular, the classical
curvature integral retains its original constant, with no
norm-equivalence factor in the conclusion. Applying this argument to
a subsequence realizing the lower limit proves the assertion for the
original sequence.
\end{proof}

\begin{remark}
The preceding lemma uses only curvature-varifold compactness.
The strong $L^2$ statement below additionally uses the vanishing
shrinker defect and the strict inequality $n>2$.
\end{remark}

\begin{lemma}\label{lem:curvature}
Let \(n>2\), let \(U\subset\R^{n+1}\) be open, and let
\(\Sigma\subset U\) be a smooth embedded self-shrinker. Suppose that
\(M_j\subset U\) are smooth, properly embedded hypersurfaces without
boundary in \(U\), and that
\[
  |M_j|\rightharpoonup \theta|\Sigma|
  \qquad\text{as varifolds in }U.
\]
Let \(\mu_j\) denote the area measure on \(M_j\), and define the
shrinker defect by
\[
  \Res_j(x)=\HH_j(x)+\tfrac12 P_j^\perp(x)x,
\]
where \(P_j^\perp(x)\) denotes orthogonal projection onto the normal
line to \(M_j\) at \(x\). Assume that for every compact set
\(K\subset U\),
\[
  \Lambda_K:=\sup_j\int_{M_j\cap K}|A_j|^n\dd\mu_j<\infty,
  \qquad
  \lim_{j\to\infty}\int_{M_j\cap K}|\Res_j|^2\dd\mu_j=0.
\]

Let \(\mathcal T\subset U\) be a tubular neighborhood of a portion of
\(\Sigma\), with nearest-point projection \(\pi\). Then for every
compact set \(K\subset\mathcal T\),
\begin{equation}\label{eq:strongB}
  \lim_{j\to\infty}
  \int_{M_j\cap K}
    \bigl|B_j(x)-B_\Sigma(\pi(x))\bigr|^2\dd\mu_j(x)=0.
\end{equation}
Here both second fundamental forms are viewed as
\(\R^{n+1}\)-valued bilinear forms on \(\R^{n+1}\) by projecting
their arguments onto \(T_xM_j\) and \(T_{\pi(x)}\Sigma\),
respectively.

For a family \(M_{j,t}\) on a fixed interval \(I=[a,b]\), suppose
that the spacetime varifolds converge to \(dt\otimes\theta|\Sigma|\) and that,
for every compact set \(K\subset U\),
\[
  \Lambda_K:=\sup_j\sup_{t\in I}
     \int_{M_{j,t}\cap K}|A_{j,t}|^n\dd\mu_{j,t}<\infty,
  \qquad
  \lim_{j\to\infty}
  \int_I\int_{M_{j,t}\cap K}|\Res_{j,t}|^2\dd\mu_{j,t}\,dt=0.
\]
Then for every compact set \(K\subset\mathcal T\),
\[
  \lim_{j\to\infty}
  \int_I\int_{M_{j,t}\cap K}
    \bigl|B_{j,t}(x)-B_\Sigma(\pi(x))\bigr|^2\dd\mu_{j,t}(x)\,dt=0.
\]
\end{lemma}

\begin{proof}
Work first in a two-sided tubular chart. Let $d$ denote signed distance
to $\Sigma$, let $N=\bar\nabla d$, and let $\pi$ be nearest-point
projection. For a compact set $K$ in the chart, choose a smooth cutoff
$0\le\chi\le1$ equal to one near $K$ and supported in a compact subset
$K'$ of the chart. Define
\[
 f_j=d|_{M_j},\qquad z_j=\nabla^{M_j}f_j=P_jN.
\]
On the limiting varifold, $x\in\Sigma$ and $P=T_x\Sigma$, so $PN=0$.
The function $(x,P)\mapsto|PN(x)|^q$ is continuous and bounded on
compact subsets of the chart for every $q>0$. Inserting a cutoff equal
to one on $K'$ therefore gives
\begin{equation}\label{eq:tilt}
 \int_{M_j\cap K'}|z_j|^q\dd\mu_j\longrightarrow0
 \qquad(q>0).
\end{equation}
All integrals until the time-dependent part of the proof are over $M_j$.

For tangent vectors $u,v$, differentiation of $f_j$ gives
\[
 \nabla^2f_j(u,v)
 =\bar\nabla^2d(u,v)+\langle B_j(u,v),N\rangle.
\]
Taking the trace and using the definition of $\Res_j$, we obtain
\begin{align}\label{eq:distance-lap}
 \Delta f_j
 &=\tr_{P_j}\bar\nabla^2d+\HH_j\cdot N\notag\\
 &=F_d(x,P_j)+\Res_j\cdot N,
\end{align}
where
\[
 F_d(x,P)=\tr_P\bar\nabla^2d
                  -\tfrac12(P^\perp x)\cdot N.
\]
On $\Sigma$, $\bar\nabla^2d|_{T\Sigma}=S_\Sigma$ and the shrinker
equation reads $H_\Sigma=\langle x,N\rangle/2$. Thus
$F_d(x,T_x\Sigma)=0$. Varifold convergence and the residual assumption
imply
\begin{equation}\label{eq:lap-small}
 \int\chi^2(\Delta f_j)^2\dd\mu_j
 \le2\int\chi^2F_d(x,P_j)^2\dd\mu_j
      +2\int\chi^2|\Res_j|^2\dd\mu_j\longrightarrow0.
\end{equation}

The Bochner identity is
\[
 \tfrac12\Delta|z_j|^2
 =|\nabla^2f_j|^2+\langle z_j,\nabla\Delta f_j\rangle
                         +\operatorname{Ric}_{M_j}(z_j,z_j).
\]
Since $\nabla|z_j|^2=2\nabla^2f_j(z_j,\cdot)$, integration against
$\chi^2$ gives
\begin{align*}
 \int\chi^2|\nabla^2f_j|^2
 ={}&-2\int\chi\nabla^2f_j(\nabla\chi,z_j)
      -\int\chi^2\langle z_j,\nabla\Delta f_j\rangle
      -\int\chi^2\operatorname{Ric}_{M_j}(z_j,z_j).
\end{align*}
The middle term satisfies
\begin{align*}
 -\int\chi^2\langle z_j,\nabla\Delta f_j\rangle
 &=\int\Delta f_j\,\Div_{M_j}(\chi^2z_j)\\
 &=\int\chi^2(\Delta f_j)^2
       +2\int\chi\Delta f_j\langle\nabla\chi,z_j\rangle.
\end{align*}
The Gauss equation and $|\HH_j|\le\sqrt n|A_j|$ give
\begin{align*}
 \operatorname{Ric}_{M_j}(z_j,z_j)
 &=\HH_j\cdot B_j(z_j,z_j)-\sum_i|B_j(z_j,e_i)|^2\\
 &\ge-(\sqrt n+1)|A_j|^2|z_j|^2.
\end{align*}
We also have
\begin{align*}
 2\chi|\nabla^2f_j||\nabla\chi||z_j|
 &\le\tfrac12\chi^2|\nabla^2f_j|^2
                         +2|\nabla\chi|^2|z_j|^2,\\
 2\chi|\Delta f_j||\nabla\chi||z_j|
 &\le\chi^2(\Delta f_j)^2+|\nabla\chi|^2|z_j|^2.
\end{align*}
Absorbing the half-Hessian term and multiplying by two yields
\begin{equation}\label{eq:bochner-est}
 \int\chi^2|\nabla^2f_j|^2
 \le4\int\chi^2(\Delta f_j)^2
       +6\int|\nabla\chi|^2|z_j|^2
       +C_n\int\chi^2|A_j|^2|z_j|^2.
\end{equation}
H\"older's inequality gives
\begin{equation}\label{eq:holder-tilt}
 \int\chi^2|A_j|^2|z_j|^2
 \le\left(\int_{M_j\cap K'}|A_j|^n\right)^{2/n}
     \left(\int_{M_j\cap K'}|z_j|^{2n/(n-2)}\right)^{(n-2)/n}
 \longrightarrow0.
\end{equation}
Equations
\eqref{eq:tilt}, \eqref{eq:lap-small}, and \eqref{eq:holder-tilt}
show that the left-hand side of \eqref{eq:bochner-est} tends to zero.

Define the ambient normal-valued bilinear tensor
\[
 C_d(x,P)=-(P\bar\nabla^2dP)\otimes N.
\]
The component in the output direction $N$ of $B_j-C_d(x,P_j)$ is
$\nabla^2f_j\otimes N$. To compute the perpendicular component, choose
either unit normal $\nu_j$ and define $b_j$ by
$B_j(u,v)=b_j(u,v)\nu_j$. Then $|b_j|=|A_j|$ and
\[
 (B_j-C_d)^{\perp N}(u,v)
 =b_j(u,v)\bigl(\nu_j-(\nu_j\cdot N)N\bigr).
\]
Since $|N|=|\nu_j|=1$,
\[
 \bigl|\nu_j-(\nu_j\cdot N)N\bigr|^2
 =1-(\nu_j\cdot N)^2=|P_jN|^2.
\]
The two output components are orthogonal. Summing over an orthonormal
basis of the input tangent space therefore gives the exact identity
\begin{equation}\label{eq:tensor-identity}
 |B_j-C_d(x,P_j)|^2
 =|\nabla^2f_j|^2+|A_j|^2|P_jN|^2.
\end{equation}
On the limiting support,
$C_d(x,T_x\Sigma)=-S_\Sigma\otimes N=B_\Sigma(x)$.
Consequently the continuous, locally bounded integrand
\[
 |C_d(x,P)-B_\Sigma(\pi x)|^2
\]
vanishes on the limit. Its integral against a compact cutoff tends to
zero. Together with \eqref{eq:tensor-identity}, this proves
\eqref{eq:strongB}. To estimate the curvature away from $\Sigma$, fix
a compact set $K\Subset U$ and an open neighborhood $\mathcal O$ of
$\Sigma\cap K$ such that $K\setminus\mathcal O$ is a positive distance
from $\Sigma$. Let
\[
 E_j:=M_j\cap(K\setminus\mathcal O).
\]
The limiting mass measure vanishes on the compact set $K\setminus\mathcal O$,
so varifold convergence gives $\mu_j(E_j)\to0$. H\"older's inequality now gives
\begin{equation}\label{eq:outside}
 \int_{E_j}|A_j|^2\dd\mu_j
 \le\Lambda_K^{2/n}\mu_j(E_j)^{1-2/n}\longrightarrow0.
\end{equation}
Finitely many tubular charts cover each compact portion of $\Sigma$.

For the spacetime assertion, use the same cutoff $\chi$ and compact set $K'$ in the tubular chart. Set
\[
 e_j(u)=\int_{M_{j,u}\cap K'}|P_jN|^{2n/(n-2)}\dd\mu_{j,u}.
\]
Spacetime varifold convergence gives $\int_Ie_j(u)\,du\to0$.
Applying H\"older's inequality first on each slice and then on $I$, with
$\beta=(n-2)/n$, gives
\begin{align*}
 \int_I\int\chi^2|A_j|^2|P_jN|^2
 &\le\Lambda_{K'}^{2/n}\int_I e_j(u)^\beta\,du\\
 &\le\Lambda_{K'}^{2/n}|I|^{1-\beta}
                         \left(\int_Ie_j(u)\,du\right)^\beta\\
 &=\Lambda_{K'}^{2/n}|I|^{2/n}
                         \left(\int_Ie_j(u)\,du\right)^{(n-2)/n}
 \longrightarrow0.
\end{align*}
The other terms in the integrated Bochner estimate tend to zero by
spacetime varifold convergence and the residual assumption. The tensor
identity then gives strong spacetime $L^2$ convergence.

\end{proof}

\subsection{The limiting self-shrinkers}

\begin{lemma}\label{lem:models}
Assume the hypotheses of Proposition~\ref{prop:conditional}. Every
sequence $s_j\to\infty$ has a subsequence for which the shifted
normalized flows on $[0,1]$ converge in local weak spacetime-varifold
topology to $du\otimes m|\Sigma|$. The support $\Sigma$ is a connected
smooth properly embedded self-shrinker, either compact or a plane, and
\begin{equation}\label{eq:model-bounds}
 m\int_\Sigma|A_\Sigma|^n\dd\mu_\Sigma\le\Lambda,
 \qquad m\int_\Sigma\rho\dd\mu_\Sigma=\Theta.
\end{equation}
For the same subsequence, the slice varifolds converge at every sequence
of times $u_j\in[0,1]$. Moreover, the convergence is uniform in $u$ against every fixed compactly
supported continuous Grassmann-bundle test.
\end{lemma}

\begin{proof}
Let
\[
 D_j=\int_{s_j}^{s_j+1}\int_{M_s}|\Res|^2\rho\dd\mu_s\,ds.
\]
By \eqref{eq:defect}, $D_j\to0$. For a compactly supported smooth
ambient function $\phi$ and $0\le u\le v\le1$, \eqref{eq:local-mass}
and Cauchy--Schwarz give
\begin{align*}
 &\left|\int_{M_{s_j+v}}\phi\rho\dd\mu_{s_j+v}
              -\int_{M_{s_j+u}}\phi\rho\dd\mu_{s_j+u}\right|\\
 &\quad\le
 \left(\int_{s_j+u}^{s_j+v}\int
                    |\bar\nabla\phi|^2\rho\dd\mu_s\,ds\right)^{1/2}
 \left(\int_{s_j+u}^{s_j+v}\int|\Res|^2\rho\dd\mu_s\,ds\right)^{1/2}
       +\|\phi\|_\infty D_j\\
 &\quad\le E^{1/2}\|\bar\nabla\phi\|_\infty D_j^{1/2}
                         +\|\phi\|_\infty D_j.
\end{align*}
Thus
\begin{equation}\label{eq:mass-osc}
 \sup_{u,v\in[0,1]}
 \left|\int\phi\rho\dd\mu_{s_j+u}
                     -\int\phi\rho\dd\mu_{s_j+v}\right|
 \le C_\phi D_j^{1/2}+\|\phi\|_\infty D_j\longrightarrow0.
\end{equation}
Replacing $\phi$ by $\phi/\rho$ gives the same vanishing oscillation
for unweighted compactly supported mass tests.

The slice varifolds have uniformly bounded local first variation. Indeed,
\begin{equation}\label{eq:variation}
 \int_{M_s\cap B_R}|\HH|\dd\mu_s
 \le\sqrt n\left(\int_{M_s\cap B_R}|A|^n\dd\mu_s\right)^{1/n}
                 \mu_s(B_R)^{1-1/n}
 \le\sqrt n\Lambda^{1/n}(DR^n)^{1-1/n}.
\end{equation}
For a smooth vector field $Y$ supported in $B_R$,
$\delta|M_s|(Y)=-\int\HH\cdot Y\dd\mu_s$.
Thus \eqref{eq:area}, \eqref{eq:variation}, and Allard's compactness
theorem \cite{Allard} apply to arbitrary sequences of slices.

Put $\lambda_j=e^{s_j/2}$ and consider
\[
 \mathcal N_j(\tau)
 =\lambda_j\bigl(M_{T+\lambda_j^{-2}\tau}-x_0\bigr),
 \qquad -\lambda_j^2T<\tau<0.
\]
Ilmanen's compactness theorem \cite[Theorem~7.1]{IlmanenCompact} and a diagonal subsequence on compact negative-time
intervals give a tangent flow. Huisken's monotonicity formula makes this
flow backward self-similar. The uniform Gaussian tail bound above and
local convergence at a negative time of convergence identify the
Gaussian mass of its normalized slice with \(\Theta>0\). Thus the tangent
flow is nonzero. Its support at negative times is smooth and
properly embedded by the hypothesis of the proposition. Under
$\tau=-e^{-u}$ we have
\[
 (-\tau)^{-1/2}\mathcal N_j(\tau)
 =e^{u/2}\lambda_j\bigl(M_{T-e^{-s_j-u}}-x_0\bigr)
 =M_{s_j+u}.
\]
The normalized limiting mass measure is therefore independent of $u$;
denote it by $\mu_\infty$. Initially the convergence of masses is known
for almost every $u$. Fix one such time. Equation \eqref{eq:mass-osc}
then gives convergence to $\mu_\infty$ at every sequence
$u_j\in[0,1]$.

For completeness, this also identifies tangent planes at arbitrary
times. Every sequence $|M_{s_j+u_j}|$ has an integral-varifold
subsequence by the preceding first-variation bound, and its mass is
$\mu_\infty$. An integral varifold supported on a smooth embedded
hypersurface has the tangent plane of that hypersurface almost
everywhere. Hence the mass determines the varifold, and every such
subsequence has the same limit $V$. It follows that
\[
 |M_{s_j+u_j}|\rightharpoonup V
 \quad\hbox{for every choice of }u_j\in[0,1].
\]
In fact, for every compactly supported continuous Grassmann-bundle
test $\Psi$,
\[
 \sup_{u\in[0,1]}
 \left|\int\Psi\dd|M_{s_j+u}|-\int\Psi\dd V\right|
 \longrightarrow0;
\]
otherwise a sequence of violating times would contradict the preceding
statement. Integrating in $u$, and approximating time-dependent tests by
finite sums of products, proves convergence to $du\otimes V$.

Let $\Sigma=\operatorname{spt}\mu_\infty$. Then
$V=\theta|\Sigma|$, where $\theta$ denotes the multiplicity of $V$. The support is closed as a subset of the ambient
space. Since it is a smooth embedded hypersurface, its inclusion is
proper: its intersection with every closed ambient ball is compact in
the submanifold topology. Each connected component is likewise proper
and complete. The integer function $\theta$ is nonnegative, and the
measure $\theta\dd\mu_\Sigma$ has support $\Sigma$.

Choose $u_j\in(0,1)$ such that
\[
 \int_{M_{s_j+u_j}}|\Res|^2\rho\dd\mu_{s_j+u_j}\le2D_j.
\]
This is possible by averaging; if $D_j=0$, any interior time may be
chosen. For every fixed $R$,
\[
 \int_{M_{s_j+u_j}\cap B_R}|\Res|^2\dd\mu_{s_j+u_j}
 \le2(4\pi)^{n/2}e^{R^2/4}D_j\longrightarrow0.
\]
For a smooth compactly supported ambient vector field $Y$,
\begin{align*}
 \int\Div_{P_j}Y\dd\mu_{s_j+u_j}
 &=-\int\HH_j\cdot Y\dd\mu_{s_j+u_j}\\
 &=\frac12\int P_j^\perp x\cdot Y\dd\mu_{s_j+u_j}
                         -\int\Res_j\cdot Y\dd\mu_{s_j+u_j}.
\end{align*}
The last term tends to zero by Cauchy--Schwarz and the local area bound.
The other integrands are continuous Grassmann-bundle functions.
Passing to the limit gives
\begin{equation}\label{eq:limit-first-variation}
 \int_\Sigma\theta\Div_\Sigma Y\dd\mu_\Sigma
 =\frac12\int_\Sigma\theta x^\perp\cdot Y\dd\mu_\Sigma.
\end{equation}
If $Y$ is tangent to $\Sigma$, the right-hand side vanishes. Every
compactly supported smooth tangential field on a coordinate patch
extends to an ambient field. Thus the distributional tangential
derivative of $\theta$ is zero. In local coordinates this says that all
weak first derivatives of the locally integrable function $\theta$
vanish; mollification shows that $\theta$ is constant on each connected
patch, and overlapping patches give constancy on each connected
component. Its value there is a positive integer, since a component
with value zero would not belong to the support.

On a component with constant multiplicity, the usual first-variation
formula changes \eqref{eq:limit-first-variation} into
\[
 -\int\HH_\Sigma\cdot Y\dd\mu_\Sigma
 =\frac12\int x^\perp\cdot Y\dd\mu_\Sigma.
\]
Arbitrary normal tests imply
$\HH_\Sigma+x^\perp/2=0$. Lemma~\ref{lem:curvature-lsc}, with $p=n$, applies to the chosen
slices. Their curvature-varifold boundary measure is zero, their local
masses satisfy \eqref{eq:area}, and their classical $L^n$ curvature
integrals satisfy \eqref{eq:total}. Taking nonnegative compact cutoffs
increasing to one gives
\[
 \int_\Sigma\theta|A_\Sigma|^n\dd\mu_\Sigma\le\Lambda.
\]
In particular, each component has finite total curvature. On every
noncompact component $\Sigma'$,
\[
 \int_{\Sigma'\cap(B_{2r}\setminus B_r)}|A_{\Sigma'}|^n\dd\mu_{\Sigma'}
 \le\int_{\Sigma'\setminus B_r}|A_{\Sigma'}|^n\dd\mu_{\Sigma'}
 \longrightarrow0.
\]
Ding's rigidity theorem \cite[Theorem~4.4]{Ding} applies to this
proper noncompact component, with empty boundary and the displayed
annular curvature decay. Its convention is
$\HH+x^\perp/2=0$, the same as ours. It follows that $\Sigma'$ is a plane
through the origin. Two distinct such planes in
$\R^{n+1}$ intersect, so they cannot be disjoint components of $\Sigma$.

There also cannot be two compact components, or a compact component
and a plane. Let $\Sigma_1,\Sigma_2$ be such a pair and set
$d_0=\operatorname{dist}(\Sigma_1,\Sigma_2)>0$.
The minimum is attained because one component is compact and the other
is closed and proper. Their homothetic flows at time $\tau<0$ have
distance
\[
 d(\tau)=\sqrt{-\tau}\,d_0,
 \qquad d'(\tau)=-\frac{d_0}{2\sqrt{-\tau}}<0.
\]
Fix a minimizing pair $x,y$ at a negative time and put $y-x=d(\tau)e$.
First variation of distance gives
$T_x\Sigma_{1,\tau}=T_y\Sigma_{2,\tau}=e^\perp$.
For an orthonormal basis $v_1,\ldots,v_n$ of $e^\perp$, choose geodesics
on the two hypersurfaces with respective initial points $x,y$ and the
same initial velocity $v_i$. At the minimum of their squared distance,
\[
 0\le\frac{d^2}{da^2}\bigg|_{a=0}|y(a)-x(a)|^2
 =2d(\tau)\langle e,B_y(v_i,v_i)-B_x(v_i,v_i)\rangle.
\]
Summing in $i$ yields
$\langle e,\HH_y-\HH_x\rangle\ge0$.
Following the minimizing points by the two homothetic flows gives a
smooth upper support for $d$ at the chosen time. Its derivative equals
$d'(\tau)$, since $d$ itself is differentiable. Tangential velocities
are orthogonal to $e$, and hence
\[
 d'(\tau)=\langle e,\HH_y-\HH_x\rangle\ge0,
\]
a contradiction. Thus $\Sigma$ is connected, and $\theta=m$ is a
single positive integer.

The bound \eqref{eq:area} and the resulting Gaussian tail bound pass to the limit. Local mass convergence and then $R\to\infty$ yield
\[
 m\int_\Sigma\rho\dd\mu_\Sigma
 =\lim_j\Fcal(s_j)=\Theta.
\]
This proves \eqref{eq:model-bounds}.

\end{proof}

In the setting of Theorem~\ref{thm:main}, the smooth-support hypothesis
used above is provided by \cite[Theorem~1.1]{HW}, whose hypotheses are
smooth closed embedded hypersurfaces and a uniform global $L^n$ bound
for the second fundamental form. The local area-ratio
and curvature hypotheses in \cite[Theorem~5.1]{HW} follow here from
\eqref{eq:area} and \eqref{eq:total}. That theorem supplies smooth
embedded support. The multiplicity and the compact-or-planar conclusion
are established separately above.

\begin{lemma}\label{lem:all-multiple}
If one tangent flow at $(x_0,T)$ has multiplicity greater than one, then
every stationary normalized limit $m|\Sigma|$ furnished by
Lemma~\ref{lem:models} has $m\ge2$.
\end{lemma}

\begin{proof}
Suppose a normalized limit has $m=1$. By
Lemma~\ref{lem:models}, it generates a tangent flow with either a compact
smooth embedded shrinker or a plane as support. In the compact case,
Schulze's theorem \cite[Corollary~1.2]{Schulze} makes this tangent flow
unique. Its bounded-area-ratio hypothesis follows from
\eqref{eq:area}. Thus no multiple tangent flow can occur. In the plane
case, \eqref{eq:model-bounds} gives $\Theta=1$. White's local regularity
theorem \cite{White05} makes $(x_0,T)$ a regular point, whose tangent
flow is a multiplicity-one plane. This is again a contradiction.
\end{proof}

\begin{corollary}\label{lem:window-errors}
Let $s_j\to\infty$ and pass to a subsequence as in
Lemma~\ref{lem:models}, with limit $m|\Sigma|$. Fix $L<\infty$ and a
smooth tubular neighborhood $\mathcal T$ of the relevant compact part
of $\Sigma$, with nearest-point projection $\pi$. Choose a smaller tube
$\mathcal T_1$ whose closure in $\overline B_L$ lies inside $\mathcal T$
and whose complement in $\overline B_L$ stays a positive distance from
$\Sigma$. In this tube, let $d$ denote the signed distance to $\Sigma$ and let
$N=\bar\nabla d$. Then
\begin{align}\label{eq:window-errors}
 &\int_0^1\int_{M_{s_j+u}\cap B_L\cap\mathcal T_1}
       \bigl(d^2+|PN|^2+|B-B_\Sigma\circ\pi|^2\bigr)
                   \dd\mu_{s_j+u}\,du\notag\\
 &\qquad+
 \int_0^1\int_{M_{s_j+u}\cap(B_L\setminus\mathcal T_1)}
                   (1+|A|^2)\dd\mu_{s_j+u}\,du
 \longrightarrow0.
\end{align}
\end{corollary}

\begin{proof}
The compact-or-planar classification permits a two-sided tube on every
required compact portion. Choose cutoffs supported in a slightly larger
tube and equal to one on the inner tube in $\overline B_L$.
Height and tilt convergence follow from spacetime varifold convergence.
On each compact set, the unweighted defect tends to zero by
\eqref{eq:defect}, since $\rho$ has a positive lower bound there.
Lemma~\ref{lem:curvature} therefore gives the curvature term in the
first line. For each $u\in[0,1]$, let
\[
 E_{j,u}:=M_{s_j+u}\cap(B_L\setminus\mathcal T_1).
\]
Since $\overline B_L\setminus\mathcal T_1$ is compact and disjoint from
$\Sigma$, spacetime varifold convergence gives
$\int_0^1\mu_{s_j+u}(E_{j,u})\,du\to0$.
On each time slice, H\"older's inequality gives
\[
 \int_{E_{j,u}}|A|^2\dd\mu_{s_j+u}
 \le\Lambda^{2/n}\mu_{s_j+u}(E_{j,u})^{(n-2)/n}.
\]
Integrating in $u$ and applying H\"older's inequality on $[0,1]$ yields
the second line. All regions and cutoffs in this corollary are fixed after
the limiting model has been selected.
\end{proof}

\subsection{Connected components in ambient balls}\label{subsec:branch-detection}\leavevmode\par\nobreak

For a smooth closed embedded hypersurface \(S\) and \(z\in S\), define
\begin{equation}\label{eq:sep-def}
 \mathfrak s_S(z)=\inf\{r>0:S\cap B_r(z)\text{ is disconnected}\},
 \qquad \inf\varnothing=+\infty .
\end{equation}
At ordinary time \(t\) let \(\mathfrak s(x,t)=\mathfrak s_{M_t}(x)\).
For the normalization in \eqref{eq:normalization}, let
\begin{equation}\label{eq:sep-scale}
 \widehat{\mathfrak s}(X,s)=\frac{\mathfrak s(x,t)}{\sqrt{T-t}} .
\end{equation}
Indeed the map \(y\mapsto(y-x_0)/\sqrt{T-t}\) identifies \(M_t\cap B_r(x)\) with \(M_s\cap B_{r/\sqrt{T-t}}(X)\), preserving its connected components. In particular \(\widehat{\mathfrak s}(\cdot,s)\) is exactly the separation of the normalized slice.

The strict sublevel sets of \(\mathfrak s\) are relatively open on a smooth closed flow. To see this, suppose \(\mathfrak s(x,t)<a\). Choose \(r<a\) for which the intersection is disconnected, and a point \(y\) in a component other than the one containing \(x\). For any \(r'\in(|y-x|,r)\) these two points remain in distinct components of the smaller intersection. Choose \(r'\) so that \((r')^2\) is a regular value of the smooth squared distance from \(x\) on \(M_t\). Transversality of this level persists as the source and time vary near \((x,t)\), by compactness. The sublevel sets at this fixed radius are diffeomorphic there: on a collar of the regular level one can transport it by the gradient of the defining function divided by its squared norm, and extend this transport by a cutoff off the collar. Thus nearby intersections are disconnected as well. This proves the assertion and in particular measurability of all the integrals involving \(\mathfrak s\) below.

For \(z\in S\cap B\) denote the connected component of \(S\cap B\) containing \(z\) by \(C_z(S\cap B)\). The estimate of \cite[Corollary~4.4]{HW} has the following consequence. Given \(D<\infty\) and \(0<\eta<1\), choose \(b=b(n,D,\eta)>0\) sufficiently small. If \(S\) has no relative boundary in \(B_R(p)\), has area ratios bounded as in \eqref{eq:area}, and satisfies
\(\int_{S\cap B_R(p)}|A|^n\dd\mu_S<b\), then
\begin{equation}\label{eq:branch-area-full}
 (1-\eta)\omega_n r^n\le
 \mu_S\big(C_z(S\cap B_r(z))\big)
 \le(1+\eta)\omega_n r^n
\end{equation}
whenever \(z\in S\cap B_{R/4}(p)\) and \(0<r\le R/512\).
In that corollary the curvature hypothesis is stated in terms of
\(\|A\|_{L^n}\), and the relative error is at most \(C_n\|A\|_{L^n}^{4/3}\) on using this norm as the small parameter. It suffices to take \(b^{1/n}\) below its smallness threshold and \(C_nb^{4/(3n)}\le\eta\).

\begin{lemma}\label{lem:branch-sep}
Let \(S_j\) be smooth closed embedded hypersurfaces with the bounds \eqref{eq:total} and \eqref{eq:area}, and suppose \(|S_j|\rightharpoonup m|\Sigma|\) locally, where \(\Sigma\) is smooth and embedded and \(m\ge2\). After passage to a subsequence, there is a finite set \(Q\subset\R^{n+1}\) with the following property. For each \(p\in\Sigma\setminus Q\) and any \(z_j\in S_j\) tending to \(p\),
\[
                    \mathfrak s_{S_j}(z_j)\longrightarrow0 .
\]
\end{lemma}
\begin{proof}
The measures \(|A_j|^n\dd\mu_j\) have total mass at most \(\Lambda\). Pass to a locally weakly convergent subsequence with limit \(\sigma\), so \(\sigma(\R^{n+1})\le\Lambda\). Use \(\eta=1/4\) in \eqref{eq:branch-area-full}, and set
\[
 Q=\{p:\sigma(\{p\})\ge b\},\qquad b\,\# Q\le\Lambda .
\]
Fix \(p\in\Sigma\setminus Q\) and \(a>0\). Continuity of \(\sigma\) on decreasing closed balls gives \(R>0\) such that \(\sigma(\overline B_R(p))<b\). By weak convergence on closed sets the curvature mass on \(B_R(p)\) is then less than \(b\) for all large \(j\). Smoothness of \(\Sigma\) at \(p\) gives
\[
        \lim_{r\downarrow0}\frac{m\mu_\Sigma(B_r(p))}{\omega_n r^n}
        =m\ge2 .
\]
Choose \(0<r<\min\{a,R/512\}\) with \(m\mu_\Sigma(B_r(p))>\frac32\omega_n r^n\). By continuity of mass on increasing open balls there exists \(r'<r\) such that the same strict inequality with right-hand side \(\frac32\omega_n r^n\) holds for \(B_{r'}(p)\).
For any \(z_j\to p\) one has \(z_j\in B_{R/4}(p)\) and \(B_{r'}(p)\subset B_r(z_j)\) for all large \(j\). Thus
\begin{align*}
 \liminf_j\mu_j(S_j\cap B_r(z_j))
 &\ge \liminf_j\mu_j(B_{r'}(p))
 \ge m\mu_\Sigma(B_{r'}(p))
   >\tfrac32\omega_n r^n,\\
 \mu_j\big(C_{z_j}(S_j\cap B_r(z_j))\big)
 &\le\tfrac54\omega_n r^n .
\end{align*}
For large \(j\) the component through \(z_j\) therefore omits part of the intersection. Formula \eqref{eq:sep-def} yields \(\mathfrak s_{S_j}(z_j)\le r<a\). Since \(a>0\) was arbitrary, this proves the lemma.
\end{proof}

\section{The separation function and the local scale}\label{sec:contact}

\subsection{A local scale on the spacetime track}\label{sec:local}\leavevmode\par\nobreak

We begin in ordinary mean curvature flow coordinates.
Let $I$ be an open interval, let $\mathcal M=\{(x,t):t\in I,\ x\in M_t\}$,
and suppose that
\begin{equation}\label{eq:loc-assumptions}
 \mu_t(B_r(y))\le D r^n,\qquad
 \int_{M_t}|A|^n\dd\mu_t\le\Lambda
 \quad (y\in\R^{n+1},\ r>0,\ t\in I).
\end{equation}
For $X=(x,t)$ and $Y=(y,u)$ set
\[
 d_P(X,Y)=\max\{|x-y|,\sqrt{|t-u|}\},\qquad
 P(X,r)=B_r(x)\times(t-r^2,t+r^2).
\]
We use a scale that records both the curvature and the connectedness
of the flow near the source. Fix a dimensional number
$0<\epsilon_L<10^{-3}$, decreasing it below if necessary.
A radius $r>0$ is admissible at $X\in\mathcal M$ if
$[t-r^2,t+r^2]\subset I$ and, whenever $0<a\le r$ and $|u-t|\le a^2$,
\begin{equation}\label{eq:loc-admissible}
 M_u\cap B_a(x)\text{ is nonempty and connected},\qquad
 \sup_{M_u\cap B_a(x)} |A|\le\epsilon_L/a .
\end{equation}
The center of the spatial ball stays at $x$ as $u$ varies. Define
\begin{equation}\label{eq:loc-envelope}
 r_*(X)=\sup\{r:r\text{ is admissible at }X\},\qquad
 \ell(X)=\inf_{Y\in\mathcal M}\{r_*(Y)+d_P(X,Y)\}.
\end{equation}
Infinite values are allowed. Recall that the separation function
$\mathfrak s(x,t)$ is the infimum of the radii $a$ for which
$M_t\cap B_a(x)$ is disconnected \cite[Definition~3.1]{BK}.

\begin{lemma}\label{lem:local-scale-basics}
At every point of the smooth track,
\begin{equation}\label{eq:loc-comparison}
 0<\ell\le r_*\le\mathfrak s,\qquad
 |\ell(X)-\ell(Y)|\le d_P(X,Y)
\end{equation}
where the scales in the difference are finite.
Under a parabolic dilation by $\lambda$ all three lengths are multiplied
by $\lambda$. For a subinterval $J=(a,b)\subset I$, superscript $J$
denotes the scales on the restricted track. At $X=(x,t)$, $t\in J$,
\begin{equation}\label{eq:loc-restrict}
 d_J(t)=\min\{\sqrt{t-a},\sqrt{b-t}\},\qquad
 r_*^J(X)=\min\{r_*(X),d_J(t)\},\quad
 \ell^J(X)=\min\{\ell(X),d_J(t)\}.
\end{equation}
\end{lemma}

\begin{proof}
Fix a point of the smooth track. Proper embeddedness and the smooth
dependence on time give a cylinder about that point on which the flow
is a single graph, with bounded curvature and bounded normal velocity
on a smaller closed cylinder. In a sufficiently small neighborhood
of each source $y$ in the smaller cylinder, squared distance to $y$
has positive definite Hessian in graph coordinates. Indeed for a
graph $z\mapsto(z,h(z))$ the Hessian is
\begin{equation}\label{eq:loc-ambient-hessian}
 2I+2Dh\otimes Dh+2(h-y_{n+1})D^2h .
\end{equation}
Take the graph coordinates over a disk containing the relevant
sublevels in its interior. Uniform bounds on this smaller chart give
strict convexity on a fixed disk about the projected source, hence
connectedness of every sufficiently small ball intersection that meets
the graph there. The point at the source moves at most
$a^2\sup|\HH|<a$ during a time interval of length $a^2$ when $a$ is
small. This also ensures nonemptiness. Thus $r_*$ has a positive
lower bound on a neighborhood of the fixed point. The other sources
in \eqref{eq:loc-envelope} have parabolic distance bounded away from
zero. Taking the infimum proves positivity of $\ell$.
Using $Y=X$ proves $\ell\le r_*$. Condition \eqref{eq:loc-admissible}
at $u=t$ gives $r_*\le\mathfrak s$. Taking the infimum in
\[
 r_*(Z)+d_P(X,Z)\le d_P(X,Y)+r_*(Z)+d_P(Y,Z)
\]
and interchanging $X,Y$ proves the Lipschitz estimate. The dilation
law follows directly from the time scale $r^2$ and the scaling law for
$|A|$ in \eqref{eq:loc-admissible}.

Restriction to $J$ imposes the additional requirement
$[t-r^2,t+r^2]\subset J$. The supremum of available radii consequently
satisfies the asserted formula for $r_*^J$. Since
$|\sqrt d-\sqrt e|\le\sqrt{|d-e|}$ for $d,e\ge0$, $d_J$ is
$1$-Lipschitz in parabolic distance. In the infimum over $Y\in\mathcal M_J$,
the contribution from $d_J(t_Y)+d_P(X,Y)$ therefore has infimum
$d_J(t)$. The contribution from $r_*(Y)+d_P(X,Y)$ has infimum at
least $\ell(X)$. Whenever $\ell(X)<d_J(t)$, every sufficiently close
competitor for $\ell(X)$ already lies in $\mathcal M_J$, since points
at times outside $J$ are at distance at least $d_J(t)$. These two
observations give the second identity in \eqref{eq:loc-restrict}.
\end{proof}

\subsection{Graphical neighborhoods and the planar model}

\begin{lemma}\label{lem:loc-graph}
There exists $a_n>0$ with the following property. If $r_*(X)\ge r$
at $X=(x,t)$, the flow in $P(X,a_n r)$ belongs to a single smooth
graph over $T_xM_t$. The graph is defined on a larger tangent disk for
each $|u-t|<(a_n r)^2$. On a smaller cylinder its derivatives have
uniform bounds after parabolic rescaling by $r^{-1}$. In particular
any Brakke limit of flows with $r_{*,j}(0,0)\ge1$ has a regular
multiplicity-one patch at $(0,0)$.
\end{lemma}

\begin{proof}
Translate the source to $(0,0)$ and scale so that $r=1$. The radius
$3/4$ is admissible and $|A|\le2\epsilon_L$ on its cylinder. For a
slice $|u|<1/16$ the set $D_u=M_u\cap B_{1/2}(0)$ is connected and
nonempty. On any unit-speed geodesic $\gamma$ of this slice, put
$f(\gamma)=|\gamma|^2$. Differentiating twice gives
\[
 (f\circ\gamma)'=2\langle\gamma,\dot\gamma\rangle,\qquad
 (f\circ\gamma)''=2+2\langle\gamma,B(\dot\gamma,\dot\gamma)\rangle .
\]
Consequently throughout $D_u$
\begin{equation}\label{eq:loc-convexity}
 \nabla^2 f\ge2(1-\epsilon_L)g_u .
\end{equation}
Properness makes each sublevel $\{f\le h\}\subset D_u$, $h<1/4$,
compact. Every critical point of $f$ is a nondegenerate minimum.
Every negative gradient trajectory stays in a compact sublevel, and
its accumulation points are critical. In fact $f$ decreases along
the trajectory with derivative $-|\nabla f|^2$, so an accumulation
point with $|\nabla f|>0$ would force a fixed decrease on infinitely
many visits to a smaller neighborhood. There are only finitely many
critical points in the compact sublevel since they are isolated. The
accumulation set of the trajectory is connected and hence consists of
a single minimum. Each minimum has an open basin of attraction by
nondegeneracy. Connectedness of $D_u$ now yields a unique minimum
$p_u$. Set $f_u=f(p_u)$.

Follow a gradient trajectory in the direction of increasing $f$.
By \eqref{eq:loc-convexity},
\[
 \frac{d|\nabla f|^2}{df}
 =\frac{2\nabla^2f(\nabla f,\nabla f)}{|\nabla f|^2}
 \ge4(1-\epsilon_L).
\]
Integration from $p_u$ yields the following bounds for the gradient
and for the length of the trajectory from $p_u$ to $z$:
\begin{align}
 |\nabla f(z)|^2&\ge4(1-\epsilon_L)(f(z)-f_u),\notag\\
 {\rm length}(p_u,z)
 &=\int_{f_u}^{f(z)}\frac{dh}{|\nabla f|}
 \le\sqrt{\frac{f(z)-f_u}{1-\epsilon_L}} .
 \label{eq:loc-trajectory}
\end{align}
In particular if $z\in B_{1/4}$, its normal line makes angle at most
$2\epsilon_L |z|/\sqrt{1-\epsilon_L}$ with the normal line at $p_u$,
by integration of the Weingarten equation along this path.

We make the reference plane independent of $u$. Let $z(u)$ be the
point moving by normal velocity from $z(0)=0$. As long as it remains
in $B_{1/2}$,
\[
 |z(u)|\le2\sqrt n\,\epsilon_L|u|,\qquad
 |p_u|\le |z(u)|.
\]
The first inequality, applied up to a possible first exit, keeps the
point in $B_{1/8}$ for all sufficiently small $|u|$. To compare the
tangent planes at different times, fix $h=1/16$ and any unit vector
$v\in T_0M_0$. At time zero follow the unit-speed geodesic with
initial velocity $v$ for length $h$. Its ambient acceleration has
norm at most $2\epsilon_L$, so its endpoint $w_v(0)$ satisfies
$|w_v(0)-hv|\le\epsilon_L h^2$. Follow this endpoint under the flow.
For $|u|\le c_nh^2$ its displacement is at most
$2\sqrt n\,\epsilon_L|u|$, by the same first-exit argument. Thus
\begin{equation}\label{eq:loc-chord}
 |w_v(u)-z(u)-hv|\le\epsilon_L h^2+4\sqrt n\,\epsilon_L|u|.
\end{equation}
Choose either unit normal $N_u$ at $p_u$. If a unit-speed path of
length $L$ in $D_u$ starts at $p_u$, the normal line at arclength
$a$ makes angle at most $2\epsilon_L a$ with $N_u$. Its unit tangent
therefore has component of magnitude at most $2\epsilon_L a$ along
$N_u$. Integration bounds the normal component of the chord by
$\epsilon_L L^2$. Apply this observation to the gradient paths from
$p_u$ to $w_v(u)$ and $z(u)$. By \eqref{eq:loc-trajectory} and the
displacement bound their lengths are at most $C_n h$. Hence
\[
 |\langle w_v(u)-z(u),N_u\rangle|\le C_n\epsilon_L h^2,\qquad
 |\langle v,N_u\rangle|
          \le C_n\epsilon_L (h+|u|/h)\le C'_n\epsilon_L h .
\]
Since $v$ was any unit vector in $T_0M_0$, this bounds the tangent
projection of $N_u$ onto the fixed reference plane. The angle
estimate after \eqref{eq:loc-trajectory} now gives an angle less than
$\pi/6$ between the fixed normal line and every normal line in
$D_u\cap B_{1/4}$, on decreasing the dimensional choice of $\epsilon_L$.

Let $\pi$ be orthogonal projection to the fixed plane $T_0M_0$.
The least singular value of $d\pi|_{T_zM_u}$ is the cosine of the
angle of the two normal lines. In particular $d\pi$ is invertible
with inverse norm at most $2$ on $D_u\cap B_{1/4}$. Starting at $p_u$,
lift the straight segment from $\pi(p_u)$ to each point in the tangent
disk $B_{1/32}(\pi(p_u))$. The length of each lift is at most $1/16$.
Since $|p_u|<1/8$, properness and the local inverse function theorem
continue the lift to the endpoint. Its smooth dependence on the
endpoint gives a smooth section of $\pi$ over this disk. For
$z\in D_u\cap B_{a_n}$ the projected gradient path to $p_u$ stays in
the same disk when $a_n$ is sufficiently small, by
\eqref{eq:loc-trajectory}. Lifting this path from $\pi(p_u)$ through
the section and lifting it along the original gradient path give the
same lift by local uniqueness. Thus $z$ lies in the graph of the
section. This also excludes any additional sheet in the smaller ball.
Shrinking $a_n$ once more ensures a common tangent disk for every
$|u|<a_n^2$. The spatial inverse function theorem with time as a
parameter gives smooth dependence on $u$. On interior cylinders the
graph equation is uniformly parabolic, and the curvature and slope
bounds yield higher derivative bounds by graphical interior estimates
\cite{EH}.
Smooth compactness on the smaller disk now gives the last assertion
with multiplicity one.
\end{proof}

Fix $0<\varepsilon<10^{-2}$ and an integer $k\ge3$. At $X$ with
$\ell(X)<\infty$, translate $X$ to $(0,0)$ and dilate parabolically by
$\ell(X)^{-1}$. We include $X$ in $G_\varepsilon$ if a connected
sheet through $(0,0)$ is, in these coordinates, a normal graph of
parabolic $C^k$ norm less than $\varepsilon$ over
\begin{equation}\label{eq:loc-planar-domain}
 B^P_{\varepsilon^{-1}}(0)\times(-1+\varepsilon,\varepsilon^{-1})
\end{equation}
for an $n$-plane $P$ through $0$, and if on this graph the rescaled
lengths satisfy
\begin{equation}\label{eq:loc-planar-lengths}
       1-\varepsilon<\ell<1+\varepsilon,\qquad
       \mathfrak s<(1+\varepsilon)\ell .
\end{equation}
In the graph norm, a time derivative has order two. The time interval
in \eqref{eq:loc-planar-domain} must be contained in the interval of
definition of the rescaled flow. We fix $0<c_0<10^{-3}$ and, on flows
with finite scales, put
\begin{equation}\label{eq:loc-open-set}
 \mathcal U=\mathcal M\cap
       \bigcup_{X\in\mathcal M\setminus G_\varepsilon}
                            P(X,c_0\ell(X)).
\end{equation}
Thus $\mathcal U$ is relatively open, and its definition uses the
original flow before a center for Gaussian integration is chosen.

\subsection{Convergence near a plane of multiplicity one}\leavevmode\par\nobreak

We first record the time convergence available under
\eqref{eq:loc-assumptions}. For $p=1,2$, H\"older's inequality gives
\begin{equation}\label{eq:loc-mean-bounds}
 \int_{M_t\cap B_R}|\HH|^p\dd\mu_t
 \le n^{p/2}\Lambda^{p/n}(DR^n)^{1-p/n}.
\end{equation}
The variation of a compactly supported mass test
$\phi\in C_c^1(B_R)$ is
\[
 \frac{d}{dt}\int_{M_t}\phi\dd\mu_t
 =\int_{M_t}(\langle D\phi,\HH\rangle-\phi|\HH|^2)\dd\mu_t .
\]
The two cases of \eqref{eq:loc-mean-bounds} give a uniform Lipschitz
constant for this test. If a sequence satisfying the same bounds
converges as Brakke flows, the limiting mass tests have continuous
representatives on every compact interior interval. Indeed
Arzel\`a--Ascoli gives uniformly convergent subsubsequences, and the
Brakke limit identifies all of them at almost every time. It follows
that at moving times $t_j\to t$ the masses converge to the continuous
representative at $t$. The $p=1$ estimate bounds the local first
variations of these slices. Allard compactness \cite{Allard} gives
integral-varifold subsequential limits, and the limiting mass
determines their approximate tangent planes almost everywhere. Thus
the varifolds also converge at moving times. Integrating compactly
supported tests on the Grassmann bundle then gives spacetime-varifold
convergence by dominated convergence. Gaussian tests with fixed
positive time difference are allowed by \eqref{eq:loc-assumptions}
and exhaustion by balls.

\begin{lemma}\label{lem:planes}
Suppose a sequence of flows satisfying \eqref{eq:loc-assumptions} is
defined on intervals exhausting $\R$ and converges as Brakke flows to
a static union of distinct parallel affine planes $P_0,\dots,P_q$
with positive integer multiplicities $m_i$. Suppose $m_0=1$. Let
$d=\min_{i\ne0}\dist(P_0,P_i)$ if $q\ge1$ and let $d=\infty$ if
$q=0$. The convergence is smooth and single-sheeted near $P_0$ on
compact interior spacetime sets avoiding the other planes. On compact
patches of this convergent sheet,
\begin{equation}\label{eq:loc-plane-limits}
 \begin{cases}
 \ell_j\longrightarrow d,\quad \mathfrak s_j\longrightarrow d,
                                       & d<\infty,\\
 \ell_j\longrightarrow\infty,          & d=\infty,
 \end{cases}
\end{equation}
uniformly.
\end{lemma}

\begin{proof}
Let $(z_j,u_j)$ be points on the tracks with limit $(z,u)$. At a
smooth point the Gaussian density is one. Huisken monotonicity
\cite{Huisken} therefore gives, for each fixed $a>0$ and large $j$,
\[
 1\le(4\pi a)^{-n/2}\int_{M^j_{u_j-a}}
                      e^{-|y-z_j|^2/(4a)}\dd\mu^j_{u_j-a}(y).
\]
The arbitrary-time convergence just proved applies. Integrating along
each limiting plane, and using
$|y-z|^2=|y-\pi_{P_i}z|^2+\dist(z,P_i)^2$, gives
\[
 1\le\sum_{i=0}^q m_i e^{-\dist(z,P_i)^2/(4a)} .
\]
If $z$ were outside all the planes, the right side would tend to zero
as $a\downarrow0$. Thus on a compact spacetime cylinder the supports
are eventually contained in any neighborhood of the limiting planes.

On $P_0$ the limiting density is one. On a compact subcylinder a
positive distance from the other planes, take a fixed sufficiently
small backward time difference. The contribution of the other planes
to the limiting Gaussian integrals then lies below the threshold in
White's local regularity theorem. The preceding convergence and
Gaussian tails give the same strict bound for the approximating flows
in a smaller cylinder. Monotonicity controls all smaller scales.
White's theorem \cite{White05} gives uniform interior curvature
bounds there, followed by derivative bounds. The varifold limit
identifies the smooth limit and its sheet number. A finite covering
therefore gives smooth single-sheet convergence along each prescribed
compact patch of $P_0$.

Fix $0<r<d$ and take sources $X_j=(x_j,t_j)$ converging to a point
of $P_0$ at an interior time. Suppose $\ell_j(X_j)<r$ on a
subsequence. There are competitors $Y_j=(y_j,u_j)$ satisfying
\[
          r_{*,j}(Y_j)+d_P(X_j,Y_j)<r .
\]
They lie in a compact spacetime cylinder. After extraction their
spatial limit belongs to a limiting plane and is within distance
$r<d$ of $P_0$. Hence this limit is on $P_0$. In particular the
full histories $P(Y_j,r)$ stay a fixed positive distance from every
other plane. The single-sheet assertion applies in a slightly larger
neighborhood. On graphs over $P_0$, \eqref{eq:loc-ambient-hessian}
converges to $2I$ there. All sublevels cutting out $B_b(y_j)$,
$0<b\le r$, are consequently convex in graph coordinates. The graph
can be taken over a base disk slightly larger than $r$ inside the
slab of smooth convergence, so none of these sublevels meets its
coordinate boundary. Since $|\HH_j|\to0$, for $|v-u_j|\le b^2$ the
point starting from $y_j$ moves by at most
$b^2\sup|\HH_j|<b$. Also $\sup|A_j|\le\epsilon_L/r$ there for large
$j$. These facts verify all of \eqref{eq:loc-admissible}, which
contradicts $r_{*,j}(Y_j)<r$. Thus $\liminf\ell_j(X_j)\ge d$ in
the extended sense.

If $d<\infty$, fix $r>d$ and choose a nearest plane $P_i$. The
limiting ball of radius $r$ has positive mass on $P_0$ and on
$P_i$ strictly inside the ball. By arbitrary-time mass convergence,
$M^j_{t_j}\cap B_r(x_j)$ contains points in small neighborhoods of
both planes for large $j$. All the points in this ball lie in a
union of disjoint narrow slabs about the finitely many planes, by
the preceding support estimate. The intersection is disconnected,
so $\mathfrak s_j(X_j)\le r$. Let $r\downarrow d$ and use
\eqref{eq:loc-comparison}. The same argument applies to every
countersequence from a compact patch, which proves uniformity.
\end{proof}

\subsection{A static limit with finite spacetime curvature integral}

\begin{lemma}\label{lem:flat}
Suppose flows satisfying \eqref{eq:loc-assumptions} converge to a
Brakke flow on $(-\infty,b)$ whose weight $\mu_t=\mu$ is independent
of $t$. Assume there is a single constant $C<\infty$ such that for
all $R<\infty$ and compact intervals $J\subset(-\infty,b)$,
\begin{equation}\label{eq:loc-integral-budget}
  \limsup_{j\to\infty}
   \int_J\int_{M^j_t\cap B_R}|A_j|^n\dd\mu^j_t\,dt\le C .
\end{equation}
If $\mu\ne0$, its support is a finite union of parallel affine
hyperplanes, each carrying a constant positive integer multiplicity.
If the support is a cone about $z$, it is one hyperplane through $z$.
\end{lemma}

\begin{proof}
Fix $b'<b$ and set $J_k=[b'-2k,b'-k]$. Take a countable determining
family of compactly supported mass tests. On each fixed $J_k$ these
tests converge uniformly in time to their values on $\mu$. Choose
$j(k)\uparrow\infty$ so that their first $k$ errors are at most
$1/k$ and the integral over $J_k\times B_k$ in
\eqref{eq:loc-integral-budget} is at most $C+1$. Since $|J_k|=k$,
there exists $u_k\in J_k$ such that $S_k=M^{j(k)}_{u_k}$ satisfies
\begin{equation}\label{eq:loc-vanishing-slices}
 \mu_{S_k}\rightharpoonup\mu,\qquad
 \int_{S_k\cap B_k}|A|^n\dd\mu_{S_k}\le(C+1)/k .
\end{equation}
The local mass and first-variation bounds give an integral-varifold
limit $V$ with $\|V\|=\mu$. For every fixed ball and sufficiently
large $k$, H\"older's inequality gives
\[
 \int_{S_k\cap B_R}|\HH|\dd\mu_{S_k}
 \le\sqrt n\big((C+1)/k\big)^{1/n}(DR^n)^{1-1/n}\to0 .
\]
Thus $V$ is stationary.

Fix $p\in\spt\mu$ and $s>0$. Let
$\delta_k=\|A\|_{L^n(S_k\cap B_{2048s}(p))}$. By
\eqref{eq:loc-vanishing-slices}, $\delta_k\to0$. The componentwise
area estimate \cite[Corollary~4.4]{HW}, with area-ratio bound $D$,
gives for large $k$
\begin{equation}\label{eq:loc-branch-area}
 \left|\frac{\mu_{S_k}(C_y(S_k\cap B_a(y)))}{\omega_n a^n}-1\right|
 \le C_n\delta_k^{4/3}
 \qquad(y\in S_k\cap B_{2s}(p),\ 0<a\le 2s).
\end{equation}
Here $C_y$ denotes the connected component containing $y$. List
the components $D_{k,l}$ of $S_k\cap B_s(p)$ which meet $B_{s/2}(p)$.
If $y\in D_{k,l}\cap B_{s/2}(p)$ then
$C_y(S_k\cap B_{s/4}(y))\subset D_{k,l}$. Thus each listed
component has mass at least $\frac12\omega_n(s/4)^n$ for large $k$.
Their number is at most $2\cdot4^nD/\omega_n$. Pass to a
subsequence so that this number is fixed. Pass once more so that each
$|D_{k,l}|$ has an integral-varifold limit $V_l$ relative to $B_s(p)$.
Each component is proper without relative boundary there. Its first
variation on a compactly supported vector field is bounded by the
supremum norm of the field times $\int_{S_k\cap B_s(p)}|\HH|$, so
$V_l$ is stationary. On $B_{s/2}(p)$ their sum equals $V$.

If $p\in\spt\|V_l\|$, choose $y_k\in D_{k,l}$ tending to $p$. The
entire connected set $D_{k,l}$ lies in
$C_{y_k}(S_k\cap B_{s+|y_k-p|}(y_k))$. The upper bound in
\eqref{eq:loc-branch-area} therefore gives
\[
 \|V_l\|(B_s(p))
 \le\liminf_{k\to\infty}\mu_{S_k}(D_{k,l})
 \le \limsup_{k\to\infty}
      (1+C_n\delta_k^{4/3})\omega_n(s+|y_k-p|)^n
 =\omega_n s^n .
\]
For a stationary integral varifold, the density is at least one
at each support point, by stationary monotonicity and upper
semicontinuity of the density. Centered area ratios are nondecreasing.
Consequently every ratio $\|V_l\|(B_a(p))/(\omega_n a^n)$ for
$0<a<s$ equals one. Equality in stationary monotonicity makes $V_l$
conical about $p$ inside the ball. For clarity, its tangent cone $K$
has density one at the vertex and therefore area ratio at infinity
equal to one. At every $q\in\spt\|K\|$, the area ratio at infinity is
still one by the inclusions
$B_{R-|q|}(0)\subset B_R(q)\subset B_{R+|q|}(0)$. Monotonicity
then gives equality at every radius about $q$. Hence $K$ is also
conical about $q$, and compositions of dilations about $0$ and $q$
give translations along the line spanned by $q$. The support is
therefore its linear span, and the weight is translation invariant
on that span. A locally finite translation-invariant measure on a
linear subspace is a constant multiple of its Lebesgue measure.
Since the mass of $B_R(0)$ is $\omega_n R^n$ for every $R$, the
span has dimension $n$ and the factor is one. This proves
that $V_l$ is a multiplicity-one planar disk near $p$. The finitely
many limiting branches omitting $p$ can be excluded by shrinking
the neighborhood.

We verify that disks through $p$ from different components have the
same plane. Fix two branches $D_{k,l},D_{k,h}$. A connected proper
embedded hypersurface without relative boundary separates an open
ball into two sides. One way to see this is to count crossings
modulo two of transverse paths from a fixed point. Every transverse
closed loop bounds a disk in the ball. After a relative transverse
perturbation, the preimage of the hypersurface in the disk is a
compact one-manifold with an even number of boundary points.
Properness and the absence of relative boundary justify this parity
count even when the branch is noncompact in the ball. Thus varying
the path with endpoints fixed preserves the parity. The two local
sides along a tubular neighborhood are
therefore distinguished globally. Each of these two local sides
is connected along the hypersurface by its connectedness. Every
point of the complement can be joined to such a neighborhood by a
path stopped before its first crossing. This gives exactly two
complementary components.

For each branch choose the open side disjoint from the other
branch, and denote the choices by $E_{k,l}$ and $E_{k,h}$.
These two open sets are disjoint. Indeed a point on the chosen
side of $D_{k,l}$ can be joined to $D_{k,l}$ within that side up
to its final boundary point. The whole path avoids $D_{k,h}$,
so it belongs to the complementary side of $D_{k,h}$ containing
$D_{k,l}$. Relative to $B_s(p)$ the perimeter measure of
$E_{k,l}$ is $\mu_{S_k}\lfloor D_{k,l}$, and likewise for $h$.
Local boundedness of these perimeters gives subsequential local
$L^1$ convergence of the two characteristic functions to characteristic
functions of disjoint sets $E_l,E_h$ up to null sets.

The mod-two boundary of $[E_{k,l}]$ is $[D_{k,l}]$. Local $L^1$
convergence gives local flat convergence of the full-dimensional
chains and of their boundaries. On the other hand White's
varifold-to-flat-chain theorem \cite[Theorem~3.3]{White09} applies
to $|D_{k,l}|$: its first variations are locally bounded and its
relative mod-two boundaries vanish. It follows that
$\partial[E_l]=[V_l]$ in this ball. The same conclusion holds for
$h$. Near $p$ each $V_l$ under consideration is a single planar
disk. On compact subsets of either open half-ball on its sides the
approximating perimeters tend to zero, by the branch-varifold
convergence. Lower semicontinuity of perimeter thus makes
$\mathbf1_{E_l}$ constant on each half-ball. The two constants
must differ, since $\partial[E_l]$ is the nonzero mod-two disk.
Thus each of $E_l,E_h$ occupies one side of its respective disk.
If their planes through $p$ differed, any choices of a half-ball
bounded by the first plane and a half-ball bounded by the second
would have overlap of positive volume. This contradicts their
disjointness.

It follows at every $p\in\spt\mu$ that the support is locally a
single planar disk carrying a constant positive integer multiplicity.
Along a connected component the plane and the multiplicity are
locally constant. Since the support and its components are closed
in ambient space, each component is both open and closed in its
plane. Hence it is the whole plane. Distinct such planes are
disjoint, hence parallel. For any finite subfamily with signed
offsets $b_i$ from the origin and multiplicities $m_i$, the area
bound gives for $R>\max_i|b_i|$
\[
  \omega_n\sum_i m_i(1-b_i^2/R^2)^{n/2}
   =R^{-n}\sum_i m_i\mathcal H^n(P_i\cap B_R)
   \le R^{-n}\mu(B_R)\le D .
\]
Letting $R\to\infty$ gives $\sum_i m_i\le D/\omega_n$ for every
finite subfamily, proving finiteness. If the support is conical
about $z$, dilations about $z$ preserve its finite set of parallel
planes. Their signed offsets from $z$ must therefore all vanish,
and the support is a single plane through $z$.
\end{proof}

\subsection{Evolution on a graphical disk}\leavevmode\par\nobreak

We first work in the original flow, with the parametrization
\(\partial_tF=\HH\), and set \(\Box_0=\partial_t-\Delta_{M_t}\).
The set \(G_\varepsilon\) is defined using a fixed sufficiently small
dimensional \(\varepsilon>0\). We will only use the evolution of
\(\mathfrak s\) where a disk containing the base point can be followed
for a short time. Other sheets can be arbitrarily close to this disk.

\begin{lemma}\label{lem:sep-disk}
Let \((x_*,t_*)\in G_\varepsilon\), with \(F(p_*,t_*)=x_*\).
There are an open disk \(D\subset M\) containing \(p_*\), a smaller
neighborhood \(V\Subset D\) of \(p_*\), and a time interval \(I_*\)
containing \(t_*\), such that
\begin{equation}\label{eq:sep-fixed-labels}
 \mathfrak s(F(p,t),t)
   =\min_{q\in M\setminus D}|F(p,t)-F(q,t)|
       \qquad (p\in V,\ t\in I_*).
\end{equation}
Every minimizing label belongs to the interior of \(M\setminus D\).
In material coordinates \(\log\mathfrak s\) is locally Lipschitz and
locally semiconcave jointly in the space and time variables on
\(V\times I_*\).
\end{lemma}

\begin{proof}
Translate \((x_*,t_*)\) to the origin and parabolically dilate so
that \(\ell(x_*,t_*)=1\). Use the plane in the definition of
\(G_\varepsilon\) as the horizontal plane. At the central time take
for \(D\) the labels of the graph above the open disk of radius \(5\).
For a sufficiently short time the image of \(\overline D\) remains
a graph near \(\overline B_3^n\), its boundary is outside
\(B_4(0)\), and the part of \(F(D,t)\) meeting \(B_{5/2}(0)\)
belongs to that graph. In fact these properties hold with margins
at the central time. The graphical parametrization persists by the
inverse function theorem on a compact subdisk, and the remaining
compact annulus of labels stays outside \(B_{5/2}(0)\) by continuity.
The graph function \(f(\cdot,t)\) on \(B_3^n\) has \(C^2\)-norm at
most \(2\varepsilon\) after the time interval is decreased.

Take \(x'=F(p,t)=(z',a')\) with \(p\) close to \(p_*\) and
\(|x'|<1/4\). The squared distance to \(x'\), on this graph, is
\[
 g(z)=|z-z'|^2+(f(z,t)-a')^2 .
\]
For \(v\in\R^n\) its second derivative is
\[
 D^2g(v,v)=2|v|^2+2(Df\cdot v)^2
                 +2(f-a')D^2f(v,v)\ge |v|^2
\]
if \(\varepsilon\) is small enough. For \(0<r<2\), the sublevel
\(\{g<r^2\}\) lies in \(B_{9/4}^n\), is convex, and contains \(z'\).
It follows that \(F(D,t)\cap B_r(x')\) is connected. Also
\(F(\partial D,t)\cap B_r(x')=\varnothing\). Consequently any
point of \(F(M\setminus D,t)\cap B_r(x')\) is in a different
component of \(M_t\cap B_r(x')\) from \(x'\). Indeed the portions
with labels in \(D\) and in \(M\setminus\overline D\) are disjoint
relatively open sets which exhaust this intersection.

At the central point the planar test gives
\(\mathfrak s(x_*,t_*)<1+\varepsilon\). Some ball of radius
strictly smaller than \(2\) about \(x_*\) therefore has disconnected
intersection with the slice. Such a ball meets the material
complement of \(D\). Its distance from \(x_*\) to that complement
is thus less than \(2\). The distance to the same compact set of
labels is continuous in \((p,t)\), so it remains less than \(2\)
on a smaller neighborhood. Denote it momentarily by \(d\). Balls
of radius less than \(d\) see only the connected graphical portion.
For \(d<r<2\) the ball sees both portions. This proves
\eqref{eq:sep-fixed-labels} and shows that minimizers stay away from
\(\partial D\).

On a compact smaller material chart the distances to \(M\setminus D\)
are uniformly positive. This follows from compactness of the label
sets and injectivity of each \(F(\cdot,t)\). The functions
\(\log|F(p,t)-F(q,t)|\), with \(q\in M\setminus D\), thus have a
common local Lipschitz constant and a common upper bound on their
coordinate Hessians in \((p,t)\). Their minimum has the same
Lipschitz bound. Subtracting a common quadratic polynomial with
sufficiently large positive Hessian makes every member of the
family concave on a smaller convex chart. Their infimum is then
concave as well. This proves semiconcavity.
\end{proof}

The following computation is the smooth graphical case of \cite[Proposition~3.3]{BK}.

\begin{lemma}\label{lem:sep-heat}
On the neighborhoods in Lemma~\ref{lem:sep-disk},
\(\Box_0\log\mathfrak s\ge0\) in the backward viscosity sense.
\end{lemma}

\begin{proof}
We give the computation of the separation inequality in this setting,
following \cite[Proposition~3.3]{BK}. Take a smooth function
\(\phi(p,t)\) touching \(\log\mathfrak s\) from below at
\((p,t_0)\) in a neighborhood restricted to \(t\le t_0\). Put
\(x=F(p,t_0)\). Choose a minimizing label \(q\) in
\eqref{eq:sep-fixed-labels}, set \(y=F(q,t_0)\), and denote
\[
              d=|x-y|,\qquad e=(x-y)/d.
\]
The function of two labels and time
\[
 W(p',q',t)=\log|F(p',t)-F(q',t)|-\phi(p',t)
\]
is nonnegative near \((p,q,t_0)\) for \(t\le t_0\) and vanishes at
that point. Both labels vary in open sets because \(q\) is an
interior label. Differentiating in \(q'\) gives
\(\langle e,w\rangle=0\) for \(w\in T_yM_{t_0}\).
Keeping \(p,q\) fixed, the derivative from the past satisfies
\(\partial_tW(p,q,t_0)\le0\), whence
\begin{equation}\label{eq:sep-time-touch}
       \partial_t\phi\ge d^{-1}\langle e,\HH_x-\HH_y\rangle .
\end{equation}

Choose orthonormal bases \(v_i,w_i\) of \(T_xM_{t_0},T_yM_{t_0}\),
respectively. At the fixed time \(t_0\), choose unit speed geodesics
in the two slices with initial data \((x,v_i)\) and \((y,w_i)\).
For their difference vector \(Z_i(\alpha)\) the Gauss formula gives
\[
 Z_i(0)=de,\qquad Z_i'(0)=v_i-w_i,\qquad
 Z_i''(0)=B_x(v_i,v_i)-B_y(w_i,w_i).
\]
For a nonzero \(Z\) one has
\[
 (\log |Z|)'=\frac{\langle Z,Z'\rangle}{|Z|^2},\qquad
 (\log |Z|)''
  =\frac{|Z'|^2+\langle Z,Z''\rangle}{|Z|^2}
   -2\frac{\langle Z,Z'\rangle^2}{|Z|^4}.
\]
The spatial minimum of \(W\) gives
\(\nabla^2\phi(v_i,v_i)\le(\log|Z_i|)''(0)\). Taking the trace and
using \(\tr B=\HH\), we obtain
\begin{align}
 \Delta\phi\le{}&
 \frac{\langle e,\HH_x-\HH_y\rangle}{d}
 +\frac1{d^2}\sum_{i=1}^n
       \big(|v_i-w_i|^2-2\langle e,v_i-w_i\rangle^2\big).
 \label{eq:sep-space-touch}
\end{align}
The first term cancels when we subtract \eqref{eq:sep-space-touch}
from \eqref{eq:sep-time-touch}.

It remains to choose the bases. Let \(\nu_x\) be either unit normal
at \(x\) and put \(\gamma=|\langle\nu_x,e\rangle|\). If the tangent
spaces coincide, take \(v_i=w_i\) for all \(i\); then \(\gamma=1\)
and the sum vanishes. Otherwise the tangent spaces have an
\((n-1)\)-dimensional intersection. Take \(v_i=w_i\) there for
\(1\le i<n\). The orthogonal complement of this intersection is a
two-plane containing the normals \(\nu_x,e\). Choose the signs of
the remaining two unit vectors so that
\(\langle v_n,w_n\rangle=\gamma\); the absolute cosine between the
two lines in this plane equals that between their normals. Since
\((w_n,e)\) is an orthonormal basis of the two-plane,
\[
 |v_n-w_n|^2=2(1-\gamma),\qquad
 \langle e,v_n-w_n\rangle^2=\langle e,v_n\rangle^2=1-\gamma^2.
\]
These identities also hold when the normals are perpendicular, in
which case \(\gamma=0\). Substitution gives
\begin{equation}\label{eq:sep-angle-source}
 \Box_0\phi\ge\frac{-2(1-\gamma)+2(1-\gamma^2)}{d^2}
              =\frac{2\gamma(1-\gamma)}{d^2}\ge0 .
\end{equation}
No differentiability or uniqueness of the minimizing label was used.
\end{proof}

\subsection{A smooth truncation at the local scale}\leavevmode\par\nobreak

We return to the original finite interval \(I=(0,T)\) and the open
set \(\mathcal U\) defined in \eqref{eq:loc-open-set} on its spacetime
track \(\mathcal M\). We construct a smooth
length \(b\) whose derivatives at \((x,t)\) are estimated in terms of
\(\ell(x,t)\). Its relative size will ensure that the separation
branch occurs only where Lemma~\ref{lem:sep-heat} applies.

\begin{lemma}\label{lem:sep-smooth-scale}
There exists a positive smooth ambient function \(b\) on a neighborhood
of the track such that
\begin{align}
 &\tfrac13\ell\le b\le4\ell,\qquad
 |D_x\log b|\le C\ell^{-1},\qquad
 |D_x^2\log b|+|\partial_t^{\rm amb}\log b|\le C\ell^{-2}
       \quad\text{on }\mathcal M, \label{eq:sep-b-derivatives}\\
 &b<\ell \quad\text{near each point of }\mathcal M\setminus G_\varepsilon,
       \qquad b>2\ell \quad\text{on }\mathcal M\setminus\mathcal U .
 \label{eq:sep-b-relative}
\end{align}
The constant is independent of the time interval on which \(b\) is
subsequently tested.
\end{lemma}

\begin{proof}
Extend \(\ell\) from \(\mathcal M\) to ambient spacetime by taking
\(\inf_{Y\in\mathcal M}(\ell(Y)+d_P(X,Y))\), and truncate this
extension from above at a constant greater than \(\sup_{\mathcal M}\ell\).
Such a constant exists since the defining time interval is finite.
Denote the extension by \(L_0\). It agrees with \(\ell\) on the track
and is \(1\)-Lipschitz for the parabolic distance. On its positive set
consider cylinders \(P(X,aL_0(X))\), with \(a>0\) chosen so that
\(100a<c_0\), where \(c_0\) is the fixed constant defining
\(\mathcal U\). The metric covering lemma gives disjoint cylinders
\(P_i=P(X_i,aL_i)\), \(L_i=L_0(X_i)\), whose fivefold enlargements
cover the positive set. The notation \(kP_i\) multiplies the spatial
radius by \(k\) and the time radius by \(k^2\).

For \(X\in 10P_i\), the Lipschitz bound gives
\((1-10a)L_i\le L_0(X)\le(1+10a)L_i\).
Thus at any point of overlap of tenfold cylinders the corresponding
radii are comparable. Their original disjoint cylinders lie in a
single cylinder of radius comparable to these radii. The ambient
volume of a parabolic cylinder of radius \(r\) is \(2\omega_{n+1}r^{n+3}\).
Division by the volume of a smallest such cylinder bounds the
multiplicity of the overlap. The same argument shows local finiteness
at points where \(L_0>0\).

Take nonnegative smooth cutoffs \(\zeta_i\le1\) equal to one on
\(5P_i\) and supported in \(8P_i\), obtained from a fixed bump by
the parabolic change of variables about \(X_i\) with spatial scale
\(aL_i\). Put \(\psi_i=\zeta_i/\sum_j\zeta_j\); the denominator
lies between one and the overlap bound. Differentiating this quotient
and using comparability of the radii gives
\[
 |D_x\psi_i|\le C L_i^{-1},\qquad
 |D_x^2\psi_i|+|\partial_t^{\rm amb}\psi_i|\le C L_i^{-2}.
\]
Here and below constants may depend on \(a\), which is fixed.
Mark \(i\) whenever \(10P_i\) meets \(\mathcal M\setminus G_\varepsilon\).
Define
\[
 L=\sum_i L_i\psi_i,\qquad
 k=3-\tfrac52\sum_{\text{marked }i}\psi_i,\qquad b=kL .
\]
The partition identity and scale comparison yield
\(0.9L_0<L<1.1L_0\) and \(1/2\le k\le3\). On a neighborhood of
each bad point, every nonvanishing cutoff belongs to a marked
index. Indeed each unmarked support lies in \(8P_i\) while \(10P_i\)
misses that point, and there are only finitely many supports to
consider locally. Therefore \(k=1/2\) there. Conversely, if
\(X\in\mathcal M\cap\spt\psi_i\) for a marked \(i\), choose
\(Y\in 10P_i\cap(\mathcal M\setminus G_\varepsilon)\). Then
\[
   d_P(X,Y)<18aL_i<c_0(1-10a)L_i\le c_0\ell(Y).
\]
Thus \(X\in\mathcal U\). Off \(\mathcal U\) one has \(k=3\).
These observations give the size and strict inequalities in
\eqref{eq:sep-b-derivatives}--\eqref{eq:sep-b-relative}.

At an arbitrary point the derivative bounds of the partition give
\begin{align*}
 |D_xL|&\le C,& |D_x^2L|+|\partial_t^{\rm amb}L|&\le C/L_0,\\
 |D_x k|&\le C/L_0,& |D_x^2k|+|\partial_t^{\rm amb}k|&\le C/L_0^2 .
\end{align*}
For example, each term \(L_iD_x^2\psi_i\) in \(D_x^2L\) is
bounded by \(C/L_i\le C'/L_0\), and there are at most a fixed
number of terms. Use
\(D_x^2\log L=L^{-1}D_x^2L-L^{-2}D_xL\otimes D_xL\) and
\(\partial_t^{\rm amb}\log L=L^{-1}\partial_t^{\rm amb}L\).
The corresponding identities for \(k\), together with
\(\log b=\log L+\log k\), give the remaining bounds.
\end{proof}

On the ordinary flow, for a smooth ambient field \(f=f(x,t)\) we
have
\[
 \partial_t(f\circ F)=\partial_t^{\rm amb}f+\langle D_x f,\HH\rangle,
 \qquad
 \Delta_{M_t}(f\circ F)=\tr_{TM_t} D_x^2f+\langle D_x f,\HH\rangle .
\]
Consequently \eqref{eq:sep-b-derivatives} gives
\begin{equation}\label{eq:sep-cap-heat}
 |\Box_0\log b|
  =|\partial_t^{\rm amb}\log b-\tr_{TM_t}D_x^2\log b|
  \le C\ell^{-2}.
\end{equation}

\subsection{The distributional inequality on the rescaled flow}\leavevmode\par\nobreak

Let \(\tau=T-t\), \(s=-\log\tau\), and \(X=(F-x_0)/\sqrt{\tau}\).
On \(M_s\) define \(\widehat{\mathfrak s}=\mathfrak s/\sqrt{\tau}\),
\(\widehat b=b/\sqrt{\tau}\) and \(\widehat\ell=\ell/\sqrt{\tau}\),
all evaluated at the corresponding labels and times. Denote by
\(\widehat{\mathcal U}\) the copy of \(\mathcal U\) in the normalized
flow. For \(0<\delta<1\) define
\begin{equation}\label{eq:sep-truncation}
   h_\delta=\log\min\{\delta,\widehat{\mathfrak s},\widehat b\}.
\end{equation}
This is finite even on a slice where \(\mathfrak s\) takes the
value \(+\infty\).

\begin{lemma}\label{lem:sep-weak}
In material coordinates for the normalized flow, \(h_\delta\) is
locally Lipschitz and locally semiconcave in spacetime. With
\(\Box=\partial_s-\Delta_{M_s}\), it satisfies
\begin{equation}\label{eq:weak-separation-source}
 \Box h_\delta\ \ge\
  \tfrac12\mathbf1_{\{\widehat{\mathfrak s}<\delta\}
                                \setminus\widehat{\mathcal U}}
  -C\widehat\ell^{-2}\mathbf1_{\widehat{\mathcal U}
                               \cap\{\widehat\ell<4\delta\}}
\end{equation}
in distributions. The constant is independent of \(\delta\).
\end{lemma}

\begin{proof}
On a neighborhood of each point of \(\mathcal M\setminus G_\varepsilon\)
the strict inequality \(b<\ell\le\mathfrak s\) makes
\eqref{eq:sep-truncation} the minimum of two smooth functions. In a
neighborhood of each point of \(G_\varepsilon\), Lemma~\ref{lem:sep-disk}
gives a compact family of smooth logarithmic distances with uniformly
bounded first and second derivatives on smaller charts. The smooth
change of time variable and the addition of \(s/2\) preserve these
local bounds. The minimum with the other two smooth branches is
therefore locally Lipschitz and locally semiconcave throughout the
normalized spacetime track.

For a scalar field pulled back from the original flow,
\(\partial_s=\tau\partial_t\) along the material parametrization and
\(\Delta_{M_s}=\tau\Delta_{M_t}\), since the metric on each normalized
slice is \(\tau^{-1}\) times the original metric. It follows for a
positive length \(d\) that
\begin{equation}\label{eq:sep-log-rescaled}
 \Box\log(d/\sqrt\tau)=\tau\Box_0\log d+\tfrac12 .
\end{equation}
The same formula applies to lower touching tests with the smooth
time change. Hence on the graphical neighborhoods
\(\Box\log\widehat{\mathfrak s}\ge1/2\) in backward viscosity
sense, while \eqref{eq:sep-cap-heat} implies
\(\Box\log\widehat b\ge1/2-C\widehat\ell^{-2}\).

Suppose \(\phi\) touches \(h_\delta\) from below in the past at a
given point. It touches each branch attaining the minimum there.
Off \(\widehat{\mathcal U}\), the planar test and
\eqref{eq:sep-b-relative} give
\(\mathfrak s<(1+\varepsilon)\ell< b\). There the strict
inequality \(\widehat{\mathfrak s}<\delta\) gives
\(\Box\phi\ge1/2\); if \(\widehat{\mathfrak s}\ge\delta\) the
constant branch gives \(\Box\phi\ge0\). Inside
\(\widehat{\mathcal U}\), either the constant or separation branch
gives a nonnegative lower bound whenever it attains the minimum.
If the smooth scale branch attains the minimum, then
\(\widehat\ell/3\le\widehat b\le\delta\), so
\(\widehat\ell\le3\delta<4\delta\). In this case
\(\Box\phi\ge-C\widehat\ell^{-2}\). These bounds also cover all
ties. The source in \eqref{eq:weak-separation-source} is measurable:
off the relatively open set \(\widehat{\mathcal U}\) the separation
has the local continuous representation in Lemma~\ref{lem:sep-disk}.
It is locally bounded since \(\ell>0\) is continuous.

To pass to distributions, use local coordinates \(z=(p,s)\) on the
track. By Alexandrov's theorem for semiconcave functions,
\(h_\delta\) has a quadratic Taylor expansion at almost every \(z\).
At such a point its Taylor polynomial minus \(\eta|\cdot-z|^2\)
touches it from below on a sufficiently small neighborhood for
every \(\eta>0\). Applying the preceding bounds to this polynomial
and letting \(\eta\downarrow0\) proves \eqref{eq:weak-separation-source}
for the absolutely continuous part of the operator. In these
coordinates \(D^2h_\delta\) is a matrix of measures; its density
is the Alexandrov Hessian and its singular part is negative
semidefinite \cite{AFP}. Indeed subtracting a smooth quadratic
polynomial from \(h_\delta\) gives a concave function locally. The
singular part of its spatial Hessian is therefore negative
semidefinite as well. The coefficient matrix \(g^{ij}(p,s)\)
of the spatial Laplacian is positive definite, so
\(-g^{ij}(D_{ij}h_\delta)_{\rm sing}\) is nonnegative. The first
time derivative and the spatial first derivatives are locally
bounded functions, hence contribute no singular measures.

Equivalently, when distributions are paired with \(d\mu_s\,ds\),
every nonnegative smooth, compactly supported \(\phi\) satisfies
\[
 \langle\Box h_\delta,\phi\rangle
 =-\int h_\delta\!
   \left(\partial_s\phi+(\tfrac n2-|\HH|^2)\phi
                     +\Delta_{M_s}\phi\right)d\mu_s\,ds .
\]
Here \(\partial_s d\mu_s=(n/2-|\HH|^2)d\mu_s\) in the material
parametrization. The right side is bounded below by the integral
of \(\phi\) against the source in \eqref{eq:weak-separation-source},
as follows from the two parts of the distribution just computed.
\end{proof}

\section{Proof of multiplicity one}\label{sec:proof}

\subsection{Sheet separation and multiplicity}\label{subsec:positive-area}\leavevmode\par\nobreak

Fix the point \((x_0,T)\) in Proposition~\ref{prop:conditional}. In this section the smooth-support assumption of that proposition is in force. Suppose that at least one tangent flow has multiplicity greater than one. By Lemmas~\ref{lem:models} and~\ref{lem:all-multiple}, every subsequential normalized limit at this point has the form \(m|\Sigma|\) with \(m\ge2\) and Gaussian mass \(\Theta\).

\begin{lemma}\label{lem:late-separation}
For each fixed \(a>0\) one has
\begin{equation}\label{eq:full-active-mass}
 \lim_{s\to\infty}
 \int_{M_s\cap\{\widehat{\mathfrak s}<a\}}\rho\dd\mu_s=\Theta.
\end{equation}
In particular, if
\[
 A_s(a)=\int_s^{s+1}\int_{M_u\cap\{\widehat{\mathfrak s}<a\}}
                               \rho\dd\mu_u\,du ,
\]
then \(A_s(a)\to\Theta\).
\end{lemma}
\begin{proof}
Fix \(s_j\to\infty\) and put \(S_j=M_{s_j}\). After extraction, Lemma~\ref{lem:models} identifies the varifold limit \(m|\Sigma|\). Apply Lemma~\ref{lem:branch-sep} and define Radon measures
\[
 \beta_j={\bf1}_{\{\mathfrak s_{S_j}\ge a\}}\mu_{S_j}.
\]
They have locally bounded mass, so pass to a locally weakly convergent subsequence with limit \(\beta\). Since \(0\le\beta_j\le\mu_{S_j}\), integration against nonnegative compactly supported tests gives \(0\le\beta\le m\mu_\Sigma\).
If \(p\in\operatorname{spt}\beta\cap(\Sigma\setminus Q)\), each ball \(B_{1/k}(p)\) has positive \(\beta\)-mass. For each \(k\), the lower bound for open sets under weak convergence gives an index \(j_k>j_{k-1}\) and \(z_k\in S_{j_k}\cap B_{1/k}(p)\) where \(\mathfrak s_{S_{j_k}}(z_k)\ge a\). The proof of Lemma~\ref{lem:branch-sep} applies to this subsequence as well and gives a contradiction. Thus \(\beta\) is carried by \(Q\). This set is finite and \(m\mu_\Sigma(Q)=0\), whence \(\beta=0\).

The Gaussian tails of \(\beta_j\) are bounded by the tails of \(\mu_{S_j}\) in Section~\ref{sec:limits}. First insert a compact cutoff and then let its support exhaust \(\R^{n+1}\); this gives \(\int\rho\dd\beta_j\to0\). Every sequence of times admits such a subsequence, so the same vanishing holds along the full flow as \(s\to\infty\). Finally, by \eqref{eq:F},
\[
 \int_{\{\widehat{\mathfrak s}<a\}}\rho\dd\mu_s
   =\Fcal(s)-\int_{\{\widehat{\mathfrak s}\ge a\}}\rho\dd\mu_s
   \longrightarrow\Theta .
\]
The last assertion follows on integrating this convergence on \([s,s+1]\).
\end{proof}

\subsection{The Gaussian integral}\label{subsec:gaussian-primitive}\leavevmode\par\nobreak

Let \(h_\delta\) be the function in \eqref{eq:sep-truncation}, where \(0<\delta<1\). Define
\begin{align}
 J_\delta(s)&=\int_{M_s}h_\delta\rho\dd\mu_s,\label{eq:full-J}\\
 I_s(\kappa)&=\int_s^{s+1}\int_{M_u\cap\widehat{\mathcal U}
                    \cap\{\widehat\ell<\kappa\}}
                    \widehat\ell^{-2}\rho\dd\mu_u\,du,\qquad \kappa>0 .
 \label{eq:full-I}
\end{align}
On each compact time interval \(h_\delta\) is locally Lipschitz on the smooth track. Compactness of \(M\) gives a finite cover by such neighborhoods and a finite lower bound there. Moreover \(h_\delta\le\log\delta\), and hence
\begin{equation}\label{eq:full-Jsign}
        -\infty<J_\delta(s)\le(\log\delta)\Fcal(s)<0
\end{equation}
at each finite time. The local construction in Section~\ref{sec:contact} used the original time interval \((0,T)\). In particular \(J_\delta\) at the common endpoint of two consecutive intervals is the same quantity.

\begin{lemma}\label{lem:gaussian-increment}
There is a constant \(C\), independent of \(s,\delta\), such that
\begin{equation}\label{eq:window-increment-full}
 J_\delta(s+1)-J_\delta(s)\ge
             \tfrac12 A_s(\delta)-(C+\tfrac12\delta^2)\,I_s(4\delta).
\end{equation}
\end{lemma}

\begin{proof}
We give the distributional form of \eqref{eq:weighted-test}. The Riemannian metrics, the Gaussian weight, and the volume elements depend smoothly on time on any compact subinterval. Therefore \(J_\delta\) is locally Lipschitz. For a nonnegative smooth time test \(\eta\) supported in such an interval, the ordinary product rule for a Lipschitz function and \eqref{eq:weighted-volume} give
\begin{align*}
 -\int J_\delta(u)\eta'(u)\,du
 &=\int\!\!\int_{M_u}\eta
       (\partial_u h_\delta-h_\delta|\Res|^2)\rho\dd\mu_u\,du\\
 &\quad+\int\!\!\int_{M_u}
       \eta h_\delta\Div_{M_u}(\rho X^T/2)\dd\mu_u\,du .
\end{align*}
Integrating the divergence in space gives
\(-\int\eta\langle\nabla h_\delta,X^T/2\rangle\rho\).
On the other hand \(\nabla\rho=-\rho X^T/2\), so testing the spatial Laplacian with \(\eta\rho\) gives
\(\langle\Delta h_\delta,\eta\rho\rangle
=\int\eta\langle\nabla h_\delta,X^T/2\rangle\rho\).
Both formulas use the slice area form. Thus in the sense of distributions in time,
\[
 J_\delta'=\int_{M_u}(\Box h_\delta-h_\delta|\Res|^2)\rho\dd\mu_u ,
\]
where the first term is paired against the smooth spacetime density. The term \(-h_\delta|\Res|^2\) is nonnegative. Insert \eqref{eq:weak-separation-source} and approximate the indicator of \((s,s+1)\) by smooth time tests. The endpoint terms converge by continuity of \(J_\delta\). The lower bound on the right is locally integrable since the original flow is smooth and \(\ell>0\) has a positive minimum on every compact portion of its track.

The positive term so obtained is half the Gaussian integral of \({\bf1}_{\{\widehat{\mathfrak s}<\delta\}\setminus\widehat{\mathcal U}}\).
On \(\{\widehat{\mathfrak s}<\delta\}\) we have \(\widehat\ell\le\widehat{\mathfrak s}<\delta\). Hence
\[
 {\bf1}_{\widehat{\mathcal U}\cap\{\widehat{\mathfrak s}<\delta\}}
 \le\delta^2\widehat\ell^{-2}
        {\bf1}_{\widehat{\mathcal U}\cap\{\widehat\ell<4\delta\}} .
\]
Subtracting the integral of this expression from \(A_s(\delta)\) and retaining the negative term in \eqref{eq:weak-separation-source} proves \eqref{eq:window-increment-full}.
\end{proof}

\subsection{The source estimate at a point of unit scale}\label{sec:scale}\leavevmode\par\nobreak

We return to ordinary mean curvature flow, with time variable \(t\).
Let \(Y=(y,u)\) be a spacetime point. For \(t<u\) and \(a>0\) define
\begin{align}
 \rho_Y(x,t)&=(4\pi(u-t))^{-n/2}
                  \exp\!\left(-\frac{|x-y|^2}{4(u-t)}\right),&
 D_Y(x,t)&=\HH(x,t)+\frac{(x-y)^\perp}{2(u-t)},\notag\\
 \Theta_Y(a)&=\int_{M_{u-a}}\rho_Y(x,u-a)\dd\mu_{u-a}(x).
 \label{eq:scale-gaussian}
\end{align}
The parameter \(a\) is a backward time difference. The projection in
\(D_Y\) is taken using \(T_xM_t\). If \(0<a<b\) and the flow is
defined on the time interval in question, Huisken's formula gives
\begin{equation}\label{eq:scale-gaussian-drop}
  \Theta_Y(b)-\Theta_Y(a)
   =\int_{u-b}^{u-a}\int_{M_t}|D_Y|^2\rho_Y\dd\mu_t\,dt .
\end{equation}
We also use the weighted formula with a smooth test \(f\) on the
moving hypersurfaces. For \(\sigma=u-t\), set \(\xi=x-y\).
Using \(\partial_tF=\HH\) and \(\partial_t\dd\mu_t=-|\HH|^2\dd\mu_t\),
differentiate the weight in material coordinates to obtain
\[
 \partial_t(\rho_Y\dd\mu_t)
 =\left(\frac n{2\sigma}-\frac{\langle\xi,\HH\rangle}{2\sigma}
       -\frac{|\xi|^2}{4\sigma^2}-|\HH|^2\right)\rho_Y\dd\mu_t .
\]
As \(\Div_{M_t}\xi^T=n+\langle\xi,\HH\rangle\) and
\(\nabla^{M_t}\rho_Y=-\rho_Y\xi^T/(2\sigma)\), the density on
the right is
\(-|D_Y|^2\rho_Y+\Div_{M_t}(\rho_Y\xi^T/(2\sigma))\).
Multiply by \(f\), use the product rule and integrate the
divergence on the closed slice. A second spatial integration by parts gives
\begin{equation}\label{eq:ordinary-tested}
 \frac{d}{dt}\int_{M_t} f\rho_Y\dd\mu_t
    =\int_{M_t}\big(\Box_0 f-|D_Y|^2f\big)\rho_Y\dd\mu_t ,
 \qquad \Box_0=\partial_t-\Delta_{M_t}.
\end{equation}
The area bound \(\mu_t(B_r(x))\le Dr^n\) gives an upper bound \(E_0=E_0(n,D)\)
for every \(\Theta_Y(a)\). Indeed, the contribution of
\(k\sqrt a\le |x-y|<(k+1)\sqrt a\) is bounded by
\((4\pi)^{-n/2}D(k+1)^n e^{-k^2/4}\). This also gives uniform tails
in every passage to a limit of Gaussian integrals below.

For an integral over a spacetime subset \(Q\), the notation
\(\int_Q f\dd\mu_t\,dt\) means integration first over \(M_t\cap Q_t\)
and then over time. Fix the tolerance \(\varepsilon\) in the definition
of \(G_\varepsilon\). We establish an alternative at a point where
the local scale has been normalized to one. Its use at arbitrary scales
will require keeping track of two different centers in \eqref{eq:scale-gaussian-drop}.

\begin{lemma}\label{lem:two-charges}
For every \(K<\infty\) there is \(0<\eta<1/10\), depending only on
\(n,D,\Lambda,K,\varepsilon\) and the fixed constants in the definition of
\(\ell\), with the following property. Suppose \((0,0)\in\mathcal M\),
\(\ell(0,0)=1\), and the flow is defined on \((-T_-,t_0)\), where
\[
 T_-\ge t_0/2,\qquad t_0>\eta^{-2},\qquad |z|\le K\sqrt{t_0}.
\]
Put \(Z_0=(z,t_0)\) and
\[
 Q_\eta=B_{\eta^{-1}}(0)\times
    \big(-\min\{T_-,\eta^{-2}\},\min\{t_0/2,\eta^{-2}\}\big).
\]
Let \(Z\subset Q_\eta\cap\mathcal M\) be Borel. Suppose that for every
\(Y=(y,u)\in Z\) one has \(u-\eta^{-2}>-T_-\) and
\begin{equation}\label{eq:scale-local-density-pinching}
                  \Theta_Y(\eta^{-2})-\Theta_Y(\eta^2)\le\eta .
\end{equation}
If \((0,0)\notin G_\varepsilon\), at least one of the following inequalities holds:
\begin{align}
  \int_{Q_\eta\setminus Z}|D_{Z_0}|^2\rho_{Z_0}\dd\mu_t\,dt
                    &>\eta\,t_0^{-n/2},\label{eq:scale-local-drop}\\
  \int_{Q_\eta}|A|^n\dd\mu_t\,dt &>1.\label{eq:scale-local-curvature}
\end{align}
\end{lemma}

\begin{proof}
Suppose there were counterexamples with \(\eta_j\downarrow0\).
Then \(t_{0,j}>\eta_j^{-2}\) and \(T_{-,j}\ge t_{0,j}/2\), so the
cylinders \(Q_j=Q_{\eta_j}\) eventually contain any given compact
spacetime set. The hypotheses give
\begin{align}
 \int_{Q_j}|A_j|^n\dd\mu_{j,t}\,dt &\le1,\notag\\
 t_{0,j}^{n/2}\int_{Q_j\setminus Z_j}
        |D_{Z_{0,j}}|^2\rho_{Z_{0,j}}\dd\mu_{j,t}\,dt&\le\eta_j .
 \label{eq:scale-counter-bounds}
\end{align}
Use Brakke compactness \cite{IlmanenCompact} and pass to a subsequence
on compact time intervals. The uniform slice bounds give convergence
of mass measures at arbitrary converging times, as in the proof of
Lemma~\ref{lem:planes}. In particular the limiting mass tests are
continuous in time. The slices also converge as varifolds: the
local first variations are bounded and the weight of an integral
varifold determines its tangent plane almost everywhere. Denote the
limiting slice by \(V_t\) and its weight by \(\mu_t^\infty\).
Since \(r_{*,j}(0,0)\ge\ell_j(0,0)=1\), admissibility gives a uniform
single graphical neighborhood of the origin in spacetime and a
curvature bound there. On a smaller cylinder the convergence is
smooth with multiplicity one. Thus the limit is nonzero.

We first pass a mean curvature estimate to the limit. Put
\(q=n/(n-1)\). For a smooth compactly supported ambient vector field
\(W=W(x,t)\), first variation on each slice gives
\[
 \int\!\!\int W\cdot\HH_j\dd\mu_{j,t}\,dt
 =-\int\!\!\int \tr_P(D_xW)\,d|M_{j,t}|(x,P)\,dt .
\]
Its support is contained in \(Q_j\) for large \(j\). Since
\(|\HH_j|\le\sqrt n\,|A_j|\), the absolute value of the left side is at
most \(\sqrt n\,(\int\!\!\int |W|^q\dd\mu_{j,t}\,dt)^{1/q}\).
Pass to the limit on the right by spacetime varifold convergence.
Duality gives an \(L^n(d\mu_t^\infty\,dt)\) vector field
\(\HH_\infty\) representing the negative of the first variation.
Using fields \(\alpha(t)W(x)\) identifies it with the generalized
mean curvature on almost every slice. The bound is uniform over
compact supports, and hence
\begin{equation}\label{eq:scale-counter-H}
       \int_{\R}\int|\HH_\infty|^n\dd\mu_t^\infty\,dt
                         \le(\sqrt n)^n .
\end{equation}

\emph{Case 1. No subsequence of points \(Y_j\in Z_j\) has a finite spacetime limit.}
Fix a compact cylinder \(Q'\). Eventually \(Q'\subset Q_j\setminus Z_j\).
For \((x,t)\in Q'\), \(t_{0,j}-t\) is comparable to \(t_{0,j}\), while
\(|x-z_j|\le K\sqrt{t_{0,j}}+C(Q')\). In consequence
\[
 t_{0,j}^{n/2}\rho_{Z_{0,j}}(x,t)\ge c(n,K)>0,\qquad
 \sup_{Q'}\frac{|x-z_j|}{2(t_{0,j}-t)}\longrightarrow0 .
\]
By \eqref{eq:scale-counter-bounds} and the local spacetime mass bound,
\begin{align*}
 \int_{Q'}|\HH_j|^2\dd\mu_{j,t}\,dt
 &\le2\int_{Q'} |D_{Z_{0,j}}|^2\dd\mu_{j,t}\,dt
 +2\int_{Q'}\frac{|x-z_j|^2}{4(t_{0,j}-t)^2}\dd\mu_{j,t}\,dt
 \longrightarrow0 .
\end{align*}
For \(\phi\in C_c^1(\R^{n+1})\), and \(t_1<t_2\) in \(Q'\)'s
time interval with \(\spt\phi\) in its spatial ball, the mass change
on \([t_1,t_2]\) is the integral of
\(\bar\nabla\phi\cdot\HH_j-\phi|\HH_j|^2\). Its absolute value is at most
\[
 \|\bar\nabla\phi\|_\infty
       \left(\int_{Q'}1\dd\mu_{j,t}\,dt\right)^{1/2}
       \left(\int_{Q'}|\HH_j|^2\dd\mu_{j,t}\,dt\right)^{1/2}
 +\|\phi\|_\infty\int_{Q'}|\HH_j|^2\dd\mu_{j,t}\,dt .
\]
It tends to zero. Thus \(\mu_t^\infty\) is constant in time.
The first inequality in \eqref{eq:scale-counter-bounds} applies on
every fixed compact cylinder in the infinite past. Lemma~\ref{lem:flat}
therefore identifies the limit as a finite union of parallel
affine planes with constant integer multiplicities.

The plane through the origin has multiplicity one by the central
graphical neighborhood. If it is the only plane,
Lemma~\ref{lem:planes} gives \(\ell_j(0,0)\to\infty\), contradicting
normalization. If another plane is present, let \(d>0\) be the
distance to a nearest distinct plane. The same lemma gives
\(\ell_j(0,0)\to d\), so \(d=1\). On compact patches of the
central sheet both \(\ell_j\) and \(\mathfrak s_j\) tend uniformly
to one. Smooth single-sheet convergence holds there as well.
Apply these conclusions on a slightly larger disk and time interval
than the fixed planar test window. For large \(j\) the central sheet
is a graph with the required small norm, and satisfies
\(1-\varepsilon<\ell_j<1+\varepsilon\) and
\(\mathfrak s_j<(1+\varepsilon)\ell_j\) throughout the window.
This contradicts \((0,0)\notin G_\varepsilon\).

\emph{Case 2. After passing to a subsequence, there are
\(Y_j=(y_j,u_j)\in Z_j\) with \(Y_j\to Y=(y,u)\).}
For \(0<a<b<\infty\) fixed, one has
\(\eta_j^2<a<b<\eta_j^{-2}\) for large \(j\). Monotonicity and
\eqref{eq:scale-local-density-pinching} give
\[
 0\le\Theta_{Y_j}(b)-\Theta_{Y_j}(a)\le\eta_j .
\]
At each of the moving times \(u_j-b,u_j-a\), convergence of the
mass measures and the Gaussian tail estimate pass the corresponding
area ratios to the limit. Thus the limiting Gaussian area about \(Y\)
is constant at all backward times. It is positive because at the
smooth center \(Y_j\), \(\lim_{h\downarrow0}\Theta_{Y_j}(h)=1\);
hence \(\Theta_{Y_j}(a)\ge1\).

The equality case of Gaussian monotonicity gives a homothetically
shrinking family for \(t<u\), centered at \(y\), with
\(\HH_\infty+(x-y)^\perp/(2(u-t))=0\) almost everywhere.
To see the constancy of the rescaled slices directly, rescale about
each \(Y_j\) using \(v=-\log(u_j-t)\), and denote the normalized
flow and its measure by \(N^j_v,\nu^j_v\). Put
\(\Phi_j=\HH_{N^j_v}+p^\perp/2\) and
\(\rho(p)=(4\pi)^{-n/2}e^{-|p|^2/4}\).
On the interval \(J=[-\log b,-\log a]\), the change of scale in
\eqref{eq:scale-gaussian-drop} gives
\(\int_J\int|\Phi_j|^2\rho\dd\nu^j_v\,dv\le\eta_j\).
For any compactly supported smooth test \(\psi(p)\), the normalized
mass evolution in normal coordinates is
\[
 \frac{d}{dv}\int\psi\rho\dd\nu^j_v
 =\int\big(\bar\nabla\psi\cdot\Phi_j-\psi|\Phi_j|^2\big)\rho\dd\nu^j_v .
\]
Thus the oscillation of these tests on \(J\) is bounded by
\(\|\bar\nabla\psi\|_\infty (|J|E_0\eta_j)^{1/2}
+\|\psi\|_\infty\eta_j\). The normalized limiting mass is
independent of \(v\). The normalized first-variation identity,
integrated against compact time tests and passed to the limit
using the same \(L^2\) bound, gives the shrinker equation.

For homothetic slices, the integral of \(|\HH_\infty|^n\) is
independent of time. In fact, dilation by \(\lambda>0\) multiplies
the weight by \(\lambda^n\) and the mean curvature vector at
corresponding points by \(\lambda^{-1}\), as follows by transforming
the first-variation formula. The integral of this constant on
\((-\infty,u)\) is finite by \eqref{eq:scale-counter-H}. Therefore
the constant is zero. At almost every past time the varifold is
stationary and satisfies \((x-y)^\perp=0\) almost everywhere.
This implies dilation invariance about \(y\). Indeed, testing
stationarity with \(\psi(x)(x-y)\) gives
\[
 \int\!\left(n\psi(x)+(x-y)\cdot\bar\nabla\psi(x)\right)\dd\mu^\infty_t(x)=0 .
\]
By applying this identity to \(\psi(y+\lambda(x-y))\), the derivative
of \(\lambda^n\int\psi(y+\lambda(x-y))\dd\mu^\infty_t(x)\) is zero.
Thus a past slice is a cone. Homothety and continuity of mass tests
make the whole past flow static. Apply Lemma~\ref{lem:flat} on the
infinite past, using \eqref{eq:scale-counter-bounds}. The cone must
be a single plane \(P\) through \(y\), possibly counted more than once.

We next determine the flow at times \(t\ge u\). Choose a unit normal
\(\nu_P\) and put \(h(x)=\langle x-y,\nu_P\rangle\).
Along each smooth approximating flow one has
\[
 \partial_th=\HH_j\cdot\nu_P=\Delta h,\qquad
 (\partial_t-\Delta)h^2=-2|\nabla h|^2
                       =-2|P_{T M_{j,t}}\nu_P|^2 .
\]
Take \(t_-<u\) and \(t'>t_-\), and put \(Y'=(0,t'+1)\).
Applying \eqref{eq:ordinary-tested} with \(f=h^2\) on the closed smooth slices gives
\begin{align*}
 \left[\int h^2\rho_{Y'}\dd\mu_{j,t}\right]_{t_-}^{t'}
 &= -\int_{t_-}^{t'}\int
       \big(2|\nabla h|^2+h^2|D_{Y'}|^2\big)\rho_{Y'}\dd\mu_{j,t}\,dt
 \le0 .
\end{align*}
At both endpoints \(h^2\rho_{Y'}\) has uniform integrable tails,
by the area bound and Gaussian decay. The integral at \(t_-\)
therefore tends to zero. The limiting integral at \(t'\) is zero
as well, so \(\spt\mu^\infty_{t'}\subset P\).

For almost every \(t\), express the integral slice as
\(\theta_t\mathcal H^n\lfloor P\) with integer-valued
\(\theta_t\ge0\). In plane coordinates first variation with
\(\phi e_i\) gives
\(\int_P\theta_t\partial_i\phi=-\int_P\phi\theta_t(\HH_\infty)_i\).
Hence \(\partial_i\theta_t=\theta_t(\HH_\infty)_i\in L^1_{\rm loc}(P)\).
The weak derivative of an integer-valued \(W^{1,1}_{\rm loc}\)
function is zero. For example, composing with
\(q(v)=v-\sin(2\pi v)/(2\pi)\), which equals \(v\) and has
derivative zero at each integer, proves the assertion by the Sobolev
chain rule. Thus \(\theta_t=k(t)\) is spatially constant for almost
every \(t\). A smooth compactly supported test with integral one on
\(P\) represents \(k(t)\) by a continuous mass test. This continuous
function takes integer values almost everywhere, hence everywhere,
and is constant. All other compact mass tests agree with this
multiple of the planar measure at every time by continuity. The
central regular patch gives \(k=1\). Lemma~\ref{lem:planes} now
implies \(\ell_j(0,0)\to\infty\), contradicting \(\ell_j(0,0)=1\).
This proves the alternative.
\end{proof}

\subsection{An integral estimate for the local scale}\leavevmode\par\nobreak

We now state the estimate in the original coordinates. Its left
side is restricted to the fixed open set
\[
 \mathcal U=\mathcal M\cap
       \bigcup_{X\in\mathcal M\setminus G_\varepsilon}P(X,c_0\ell(X)).
\]
In particular the same set can be used for different Gaussian centers.
Here \(0<c_0<10^{-3}\) is the constant fixed in its definition.

\begin{theorem}\label{thm:scale}
For each \(K<\infty\) there exist \(C_K<\infty\) and
\(\kappa_K\in(0,1)\), depending only on
\(K,n,D,\Lambda,\varepsilon\) and the fixed local-scale parameters,
with the following property. Let \(Z_0=(z,t_0)\), let \(0<\tau_1<\tau_2\),
and suppose the flow is defined on \([t_0-2\tau_2,t_0)\).
Put \(\tau=t_0-t\). For \(0<\kappa\le\kappa_K\) one has
\begin{align}
 &\int_{t_0-\tau_2}^{t_0-\tau_1}\tau^{-n/2}
   \int_{\mathcal U_t\cap B_{K\sqrt\tau}(z)\cap\{\ell<\kappa\sqrt\tau\}}
           \ell^{-2}\dd\mu_t\,dt \notag\\
 &\qquad\le
 C_K\big(\Theta_{Z_0}(2\tau_2)-\Theta_{Z_0}(\tau_1/2)\big)
       +C_K\Lambda\kappa^{n-2}\log\frac{4\tau_2}{\tau_1}.
 \label{eq:scale-local-integral}
\end{align}
\end{theorem}

\begin{proof} We first construct the covering. The estimate
is invariant under parabolic dilations and translations: in its left
side the factors \(\tau^{-n/2}\) and \(\ell^{-2}\) transform inversely
to \(d\mu_t\) and \(dt\), respectively. Translate \(Z_0\) to \((0,0)\)
and dilate by \(\tau_2^{-1/2}\). We continue to denote the smaller
backward time by \(\tau_1\), so that \(0<\tau_1<1\). The target set is
\[
 \mathcal T=
 \{(x,t)\in\mathcal U:-1\le t\le-\tau_1,\quad
          |x|<K\sqrt{-t},\quad \ell(x,t)<\kappa\sqrt{-t}\}.
\]
Use Lemma~\ref{lem:two-charges} with \(2(K+1)\) in place of \(K\),
and denote its parameter by \(\eta\). Choose
\[
        0<\kappa_K\le 10^{-2}\min\{\eta,\sqrt{\varepsilon},1\}.
\]
Consider the defining cylinders \(P(X,c_0 r)\) of \(\mathcal U\)
that meet \(\mathcal T\), with \(r=\ell(X)\) and \(X=(x,t)\).
For \(Y=(y,u)\in\mathcal T\cap P(X,c_0r)\), put
\(\tau_X=-t\), \(\tau_Y=-u\). Since \(\ell\) is parabolic
\(1\)-Lipschitz and \(\ell(Y)<\kappa\sqrt{\tau_Y}\), we have
\begin{align*}
 (1-c_0)r &\le \ell(Y)<\kappa\sqrt{\tau_Y},\\
 |\tau_X-\tau_Y|&\le c_0^2r^2
      \le \frac{c_0^2\kappa^2}{(1-c_0)^2}\tau_Y .
\end{align*}
Let \(a_0=c_0^2\kappa^2/(1-c_0)^2\). Since
\(\tau_1\le\tau_Y\le1\), these inequalities and
\(|x|\le |y|+c_0r\) yield
\begin{equation}\label{eq:scale-center-position}
 (1-a_0)\tau_1\le\tau_X\le1+a_0,\qquad
 r\le2\kappa\sqrt{\tau_X},\qquad
 |x|\le2(K+1)\sqrt{\tau_X}.
\end{equation}
The radii are uniformly bounded. The metric covering lemma gives
a countable pairwise disjoint collection \(P(X_i,c_0r_i)\), with
\(X_i=(x_i,t_i)\), \(r_i=\ell(X_i)\), whose fivefold enlargements
cover \(\mathcal T\). Disjointness here is in ambient spacetime.
For each index put
\[
        \tau_i=-t_i,\qquad R_i=\eta^{-1}r_i,\qquad Q_i=P(X_i,R_i).
\]
We have \(R_i^2\le4(\kappa/\eta)^2\tau_i\le\tau_i/2500\).
Combining this bound with \eqref{eq:scale-center-position} gives
\begin{equation}\label{eq:scale-history-margin}
 Q_i\subset\R^{n+1}\times(-2,-\tau_1/2),\qquad \tau_i+5R_i^2<2.
\end{equation}

Restrict temporarily to the open interval \((-2,0)\). By the
restriction formula for the scale, its value becomes
\(\min\{\ell,d_I(t)\}\), where
\(d_I(t)=\min\{\sqrt{t+2},\sqrt{-t}\}\). At any point of the planar
test window based at \(X_i\), the change in time is no more than
\(\varepsilon^{-1}r_i^2\le 4\kappa_K^2\varepsilon^{-1}\tau_i
\le\tau_i/2500\). Since \(\tau_i\le1+a_0\), throughout this window
\(d_I\ge\frac12\sqrt{\tau_i}>(1+\varepsilon)r_i\). Thus the restriction
leaves the central scale \(r_i\) fixed. It also preserves the planar
test: the time window, the graph, and separation are unchanged,
and the only comparison of the restricted scale is with numbers
at most \((1+\varepsilon)r_i\). In particular \(X_i\notin G_\varepsilon\)
for the restricted flow.

Translate \(X_i\) to the origin and rescale by \(r_i^{-1}\).
The time interval becomes
\[
    (-T'_{-,i},t'_{0,i})
     =\left(-\frac{2-\tau_i}{r_i^2},\frac{\tau_i}{r_i^2}\right),
 \qquad z'_i=-x_i/r_i .
\]
Equation \eqref{eq:scale-center-position} gives
\(T'_{-,i}\ge t'_{0,i}/2\), \(t'_{0,i}\ge1/(4\kappa^2)>\eta^{-2}\)
and \(|z'_i|\le2(K+1)\sqrt{t'_{0,i}}\). Both endpoints are far enough
from \(0\) that \(Q_\eta\) in Lemma~\ref{lem:two-charges}
is the full cylinder of radius \(\eta^{-1}\); its pullback is \(Q_i\).
Within \(Q_i\) take
\[
 Z_i=\{Y\in Q_i\cap\mathcal M:
       \Theta_Y(R_i^2)-\Theta_Y(\eta^2 r_i^2)\le\eta\}.
\]
These quantities are defined since \eqref{eq:scale-history-margin}
places even the earliest required time strictly above \(-2\).
The map sending \((y,u)\) to a Gaussian mass at a fixed positive
backward time is continuous on this region of a smooth closed flow,
so \(Z_i\) is Borel.

If prime denotes the rescaled variables, the changes of scale are
\[
 |D'|=r_i|D|,\qquad \rho'=r_i^n\rho,\qquad
 d\mu'=r_i^{-n}d\mu,\qquad dt'=r_i^{-2}dt,\qquad |A'|=r_i|A|.
\]
Thus \(|D'|^2\rho'\,d\mu'\,dt'=|D|^2\rho\,d\mu\,dt\), while
\(|A'|^n\,d\mu'\,dt'=r_i^{-2}|A|^n\,d\mu\,dt\). For each cylinder
at least one conclusion of Lemma~\ref{lem:two-charges} holds.
Put the indices satisfying its first conclusion in \(\mathcal D\)
and the remaining indices in \(\mathcal C\). In the original
coordinates we have
\begin{align}
 E_i:=\int_{Q_i\setminus Z_i}|D_{(0,0)}|^2\rho_{(0,0)}
             \dd\mu_t\,dt
   &>\eta\,r_i^n\tau_i^{-n/2},&& i\in\mathcal D ,
       \label{eq:scale-estimate-D}\\
 F_i:=\int_{Q_i}|A|^n\dd\mu_t\,dt
   &>r_i^2,&& i\in\mathcal C .
       \label{eq:scale-estimate-C}
\end{align}

\smallskip\noindent\emph{Summation over the covering.}\enspace
On \(P(X_i,5c_0r_i)\), one has \(\ell\ge(1-5c_0)r_i\)
and \(-t\ge\tau_i/2\). The slice area is at most \(D(5c_0r_i)^n\)
and the length of the time interval is \(50c_0^2r_i^2\).
It follows that
\begin{align}
 \int_{\mathcal T\cap P(X_i,5c_0r_i)}
         (-t)^{-n/2}\ell^{-2}\dd\mu_t\,dt
 &\le (\tau_i/2)^{-n/2}
       \frac{D(5c_0r_i)^n}{(1-5c_0)^2r_i^2}\,50c_0^2r_i^2
 \notag\\
 &\le C r_i^n\tau_i^{-n/2}.
 \label{eq:scale-cover-cost}
\end{align}
We estimate the sum of the right-hand side separately on
\(\mathcal D\) and \(\mathcal C\).

Fix \(Y=(y,t)\in\mathcal M\) and write \(\tau=-t\).
If \(Y\in Q_i\), then \(|\tau_i-\tau|\le R_i^2\le\tau_i/2500\),
so \(\frac12\tau\le\tau_i\le2\tau\). Also
\(r_i\le2\kappa\sqrt{\tau_i}\le4\kappa\sqrt\tau=:L\).
Group the indices containing \(Y\) according to
\begin{equation}\label{eq:scale-dyadic}
                2^{-k-1}L<r_i\le 2^{-k}L,\qquad k=0,1,2,\ldots.
\end{equation}
For a given \(k\), the pairwise disjoint ambient cylinders
\(P(X_i,c_0r_i)\) in this group are all contained in
\(P(Y,2\eta^{-1}2^{-k}L)\), since
\(d_P(X_i,Y)<R_i\). The ordinary spacetime volume of an ambient
cylinder of radius \(a\) is \(2\omega_{n+1}a^{n+3}\).
Comparison with the volume of the smallest allowed cylinder gives
at most
\[
         N_0=\left(\frac{4}{c_0\eta}\right)^{n+3}
\]
indices in any such group. This bound applies independently of \(k\).

Suppose \(Y\in Q_i\setminus Z_i\) for an index with \eqref{eq:scale-dyadic}.
The increase of \(\Theta_Y\) on
\([\eta^2r_i^2,\eta^{-2}r_i^2]\) is greater than \(\eta\).
This interval is contained in
\[
 J_k=\left[\tfrac14\eta^2\,4^{-k}L^2,\,
                         \eta^{-2}4^{-k}L^2\right].
\]
The upper endpoint is at most \(4R_i^2\). Its full backward
history from \(Y\) is defined, by \(-t\le\tau_i+R_i^2\)
and \eqref{eq:scale-history-margin}. Choose \(m\in\mathbb N\)
with \(4^m>4\eta^{-4}\). The upper endpoint of \(J_{k+m}\)
then lies below the lower endpoint of \(J_k\). In each residue
class modulo \(m\) the intervals are disjoint, and the total
increase of the nondecreasing function \(\Theta_Y\) over them
is at most \(E_0\). Thus in each class there are at most
\(E_0/\eta\) relevant scales, each with at most \(N_0\) indices.
Consequently,
\[
       \sum_{i\in\mathcal D}{\bf1}_{Q_i\setminus Z_i}(Y)
                           \le mN_0 E_0/\eta .
\]
Insert this bound in \eqref{eq:scale-estimate-D}, and use
\eqref{eq:scale-gaussian-drop} at the original center \((0,0)\).
This yields
\begin{align}
 \sum_{i\in\mathcal D}r_i^n\tau_i^{-n/2}
 &\le\eta^{-1}\sum_{i\in\mathcal D} E_i \notag\\
 &\le C_K
      \int_{-2}^{-\tau_1/2}\int_{M_t}
                   |D_{(0,0)}|^2\rho_{(0,0)}\dd\mu_t\,dt \notag\\
 &=C_K\big(\Theta_{(0,0)}(2)-\Theta_{(0,0)}(\tau_1/2)\big).
 \label{eq:scale-drop-sum}
\end{align}

For \(i\in\mathcal C\), \eqref{eq:scale-estimate-C} gives
\(r_i^n\tau_i^{-n/2}\le r_i^{n-2}\tau_i^{-n/2}F_i\).
At a fixed \(Y\) the overlap bound at one scale gives
\begin{align*}
 \sum_{i\in\mathcal C} r_i^{n-2}\tau_i^{-n/2}{\bf1}_{Q_i}(Y)
 &\le 2^{n/2}N_0 \tau^{-n/2}
                     \sum_{k=0}^\infty(2^{-k}L)^{n-2}\\
 &= \frac{2^{n/2}N_0}{1-2^{-(n-2)}}\tau^{-n/2}
                        (4\kappa\sqrt\tau)^{n-2}
 \le C_K\kappa^{n-2}\tau^{-1}.
\end{align*}
The geometric series converges since \(n>2\). By Tonelli's theorem
and the critical slice bound,
\begin{align}
 \sum_{i\in\mathcal C} r_i^n\tau_i^{-n/2}
 &\le \int_{-2}^{-\tau_1/2}\int_{M_t}|A|^n
       \sum_{i\in\mathcal C}r_i^{n-2}\tau_i^{-n/2}{\bf1}_{Q_i}
                                     \dd\mu_t\,dt \notag\\
 &\le C_K\kappa^{n-2}
       \int_{-2}^{-\tau_1/2}\frac{1}{-t}\int_{M_t}|A|^n\dd\mu_t\,dt
 \le C_K\Lambda\kappa^{n-2}\log\frac4{\tau_1}.
 \label{eq:scale-curvature-sum-new}
\end{align}
The fivefold cylinders cover \(\mathcal T\). Summing
\eqref{eq:scale-cover-cost} and using
\eqref{eq:scale-drop-sum} and \eqref{eq:scale-curvature-sum-new}
proves the estimate with \(\tau_2=1\). Restoring \(\tau_2\) by
parabolic scaling gives \eqref{eq:scale-local-integral}.
\end{proof}

\subsection{A uniform Gaussian estimate and the conclusion}\label{subsec:uniform-conclusion}

\begin{lemma}\label{lem:global-scale-error}
For the quantities in \eqref{eq:full-I},
\begin{equation}\label{eq:uniform-scale-error}
                \lim_{\kappa\downarrow0}\limsup_{s\to\infty}
                         I_s(\kappa)=0 .
\end{equation}
\end{lemma}
\begin{proof}
Put \(\tau=T-t\), \(\tau_2=e^{-s}\) and \(\tau_1=e^{-(s+1)}\). In ordinary coordinates let
\[
  \rho_{x_0,T}(x,t)=(4\pi\tau)^{-n/2}
                 \exp\!\left(-\frac{|x-x_0|^2}{4\tau}\right).
\]
At the normalized time \(u=-\log\tau\), one has
\[
 X=\frac{x-x_0}{\sqrt\tau},\qquad \dd\mu_u=\tau^{-n/2}\dd\mu_t,
 \qquad du=\frac{dt}{\tau},\qquad
 \widehat\ell^{-2}=\frac{\tau}{\ell^2}.
\]
Multiplying the four factors of the measure and the integrand cancels the two time factors and gives
\begin{equation}\label{eq:gaussian-return}
 I_s(\kappa)=
 \int_{T-\tau_2}^{T-\tau_1}
 \int_{M_t\cap\mathcal U_t\cap\{\ell<\kappa\sqrt\tau\}}
       \ell^{-2}\rho_{x_0,T}(x,t)\dd\mu_t\,dt .
\end{equation}
For all sufficiently late \(s\), the flow is defined throughout \([T-2\tau_2,T)\), including the histories in Theorem~\ref{thm:scale}. First restrict \eqref{eq:gaussian-return} to \(B_{R\sqrt\tau}(x_0)\) for fixed \(R\ge1\). Since \(\rho_{x_0,T}\le(4\pi)^{-n/2}\tau^{-n/2}\), that theorem bounds this part by
\begin{equation}\label{eq:central-gauss}
 C_R\left[
 \Theta_{(x_0,T)}(2\tau_2)-\Theta_{(x_0,T)}(\tau_1/2)
       +\Lambda\kappa^{n-2}\log(4e)\right] ,
 \qquad 0<\kappa\le\kappa_R .
\end{equation}
Both Gaussian ratios in brackets tend to \(\Theta\). In particular their difference tends to zero for this fixed center.

We next control the part outside the ball. Let \(m\in\mathbb Z^{n+1}\),
\[
 z_m=x_0+\sqrt{\tau_1}\,m,\qquad
 C_m=z_m+\sqrt{\tau_1}[-1/2,1/2)^{n+1},\qquad
 d_n=\sqrt{n+1}/2 .
\]
These cubes partition the ambient space. On \(C_m\),
\(|x-z_m|\le d_n\sqrt{\tau_1}\) and
\(|x-x_0|\ge\sqrt{\tau_1}(|m|-d_n)_+\).
Since \(\tau_1\le\tau\le e\tau_1\), this gives
\[
 C_m\subset B_{(1+d_n)\sqrt\tau}(z_m),\qquad
 \rho_{x_0,T}(x,t)\le C_n\tau^{-n/2}e^{-c_n|m|^2}.
\]
For example, for \(|m|\ge2d_n\) take \(c_n=1/(16e)\); the finitely many remaining indices are absorbed in \(C_n\). If \(C_m\) meets the complement of \(B_{R\sqrt\tau}(x_0)\) at any time in the interval, then
\(\sqrt\tau R\le(|m|+d_n)\sqrt{\tau_1}\), so \(|m|\ge R-d_n\).

Apply Theorem~\ref{thm:scale} to each such cube with center \((z_m,T)\) and the same spatial parameter \(A_0=1+d_n\). The set \(\mathcal U\) was fixed before specifying any center. By \eqref{eq:entropy}, the drop at each translated center is at most \(E\). Consequently, for \(0<\kappa\le\min\{\kappa_{A_0},1\}\), the exterior contribution is bounded above by
\begin{equation}\label{eq:tail-gauss}
           C_{A_0}(E+\Lambda)
             \sum_{\substack{m\in\mathbb Z^{n+1}\\ |m|\ge R-d_n}}
                             e^{-c_n|m|^2}.
\end{equation}
The lattice sum is finite and its tail tends to zero as \(R\to\infty\).
Given \(\eta>0\), first choose \(R\) so that \eqref{eq:tail-gauss} is below \(\eta\). With this \(R\) fixed, choose \(\kappa\le\min\{1,\kappa_R,\kappa_{A_0}\}\) such that \(C_R\Lambda\kappa^{n-2}\log(4e)<\eta\). Finally take \(s\) large enough that the first difference in \eqref{eq:central-gauss}, multiplied by \(C_R\), is below \(\eta\). For every later \(s\) the bound is \(3\eta\); it also holds for smaller \(\kappa\), since \(I_s\) is nondecreasing in \(\kappa\). This proves \eqref{eq:uniform-scale-error}.
\end{proof}

\begin{proof}[Proof of Proposition~\ref{prop:conditional}]
Under the assumption of higher multiplicity we have \(\Theta>0\). Let \(C\) be the constant in \eqref{eq:window-increment-full}. By Lemma~\ref{lem:global-scale-error} we can fix \(0<\delta<1/4\) and a normalized time \(s_1\) so that
\[
           (C+1)I_s(4\delta)\le\Theta/8
                  \qquad\text{for every }s\ge s_1 .
\]
For this fixed \(\delta\), Lemma~\ref{lem:late-separation} gives \(s_0\ge s_1\) with \(A_s(\delta)\ge 3\Theta/4\) whenever \(s\ge s_0\). Inequality \eqref{eq:window-increment-full} therefore yields
\[
        J_\delta(s+1)-J_\delta(s)
          \ge \frac12\frac{3\Theta}{4}-\frac{\Theta}{8}
          =\frac{\Theta}{4}.
\]
At \(s_0\) the integral is finite. For each positive integer \(k\), summation on the \(k\) consecutive unit intervals starting at \(s_0\) gives
\[
             J_\delta(s_0+k)\ge J_\delta(s_0)+k\Theta/4 .
\]
For large \(k\) this contradicts \eqref{eq:full-Jsign}. The proposition follows.
\end{proof}

\begin{proof}[Proof of Theorem~\ref{thm:main}]
By \cite[Theorem~1.1]{HW}, we have the smooth-support hypothesis of
Proposition~\ref{prop:conditional} at every first-singular-time point.
Apply that proposition at each point admitting a nonzero tangent flow.
\end{proof}

\begin{proof}[Proof of Corollary~\ref{cor:unique}]
By the definition of a singular point reached by the original smooth
flow, choose $t_j\uparrow T$ and $x_j\in M_{t_j}$ with $x_j\to x_0$.
For fixed $a\in(0,T)$, we have $t_j-a>0$ for all sufficiently large
$j$. Huisken's monotonicity formula, centered at $(x_j,t_j)$, gives
\[
 (4\pi a)^{-n/2}
 \int_{M_{t_j-a}}e^{-|x-x_j|^2/(4a)}\,d\mu_{t_j-a}\ge1,
\]
because the small-scale Gaussian density at $(x_j,t_j)$ is one.
The original flow is smooth near time $T-a$.
Passing to the limit $j\to\infty$ gives the same lower bound for the Gaussian ratio centered at $(x_0,T)$.
Letting $a\downarrow0$ yields $\Theta\ge1$.
The Gaussian tail estimate obtained from \eqref{eq:area} prevents mass
from escaping in a tangent rescaling, so every tangent flow has this
same positive Gaussian mass and is nonzero.

By Theorem~\ref{thm:main} and Lemma~\ref{lem:models}, every tangent flow is generated by a multiplicity-one compact self-shrinker or a multiplicity-one plane. A plane tangent gives a regular point by \cite{White05}, so it is excluded here.
The compact tangent is unique by \cite[Corollary~1.2]{Schulze}.
Brakke compactness and uniqueness imply weak convergence of all parabolic rescalings to this tangent flow.
Local Brakke regularity \cite{White05} gives smooth multiplicity-one convergence near its smooth negative-time support.
The clearing-out lemma \cite{Brakke} excludes points of the rescaled flows on compact negative-time spacetime sets disjoint from the limiting support.
Consequently the convergence is smooth on every compact negative-time spacetime set.
Equivalently, the normalized hypersurfaces converge locally smoothly with multiplicity one to the compact self-shrinker.
\end{proof}

\end{document}